\documentclass[11pt,a4paper]{article}
\usepackage[T1]{fontenc}
\usepackage{amsmath,amssymb,amsthm,mathtools}
\usepackage{newtxtext,newtxmath}
\usepackage[a4paper,margin=28mm]{geometry}
\usepackage{enumitem,authblk,microtype}
\usepackage[numbers,sort&compress]{natbib}
\usepackage{xurl}
\usepackage{xcolor}
\usepackage[colorlinks=true,linkcolor=blue!50!black,citecolor=blue!50!black,urlcolor=blue!50!black]{hyperref}
\numberwithin{equation}{section}
\newtheorem{theorem}{Theorem}[section]
\newtheorem{lemma}[theorem]{Lemma}
\newtheorem{proposition}[theorem]{Proposition}

\theoremstyle{definition}
\newtheorem{definition}[theorem]{Definition}
\newtheorem{remark}[theorem]{Remark}
\setlist{nosep}
\title{\bfseries Existence and regularity for a non-autonomous\\ Born--Infeld equation with a confining potential}
\author[1]{Fan Yang}
\author[1]{Chengxiang Zhang\thanks{Corresponding author. E-mail: \href{mailto:zcx@bnu.edu.cn}{zcx@bnu.edu.cn}.}}
\affil[1]{School of Mathematical Sciences, Beijing Normal University, Beijing 100875, China}
\date{}

\begin{document}
\maketitle
\noindent\textbf{Abstract: }We study the non-autonomous Born--Infeld equation
\(\operatorname{div}(Du/\sqrt{1-|Du|^2})-V(x)u+|u|^{p-2}u=0\)
in \(\mathbb R^N\), where \(N\ge3\) and \(2<p<2^*=2N/(N-2)\).
For a \(C^1\) potential satisfying \(\inf V>0\), \(V(x)\to\infty\) as \(|x|\to\infty\), and
\(|x\cdot\nabla V(x)|\le C_DV(x)\) with \(0\le C_D<p-2\),
we prove the existence of a positive classical solution and infinitely many distinct nontrivial classical solutions.
The argument combines the nonsmooth monotonicity trick with compactness induced by the confining potential.
A Pohozaev-type identity yields uniform bounds for parameter-dependent critical points.
A local minimization characterization identifies these critical points with weak solutions and establishes strict spacelikeness.
Local Schauder estimates then give \(C^{2,\alpha}_{\mathrm{loc}}\) regularity.
All the solutions belong to the natural energy space \(H_V^1(\mathbb R^N)\) and satisfy \(|Du|<1\).

\tableofcontents

\section{Introduction}
This section consists of four parts. We first introduce the background and basic properties of the Born--Infeld model,
and then review nonsmooth critical point theory and the monotonicity trick.
We next formulate the non-autonomous Born--Infeld equation studied here and state the main results.
Finally, we describe the organization of the paper and fix the notation.

\subsection{Background and basic properties of the Born--Infeld model}

Classical electrodynamics is based on Maxwell's equations. In the electrostatic setting, the potential $u$ satisfies the Poisson equation
\(-\Delta u = \rho,\)
where $\rho$ is the charge density.
For \(\rho \in  L^{6/5}(\mathbb{R}^3)\), the corresponding condition \(|Du| \in L^2(\mathbb{R}^3)\) is compatible with finite energy and decay at infinity \cite{10.1063/1.1508433}.
For a point charge \(\rho = \delta_{x_0}\), the classical self-energy diverges. General $L^1$ charge densities likewise need not produce finite self-energy.

In 1934, Born and Infeld \cite{0036fa1d-f8f2-3bb6-965e-2d83cd17066e} introduced an upper bound \(b\) on the field strength and a new Lagrangian based on the principle of finite field energy. In electrostatics, this Lagrangian becomes
\[
\mathcal{L}_{\mathrm{BI}} = b^2\left(1 - \sqrt{1 - \frac{|\mathbf{E}|^2}{b^2}}\right),
\]
After variation and normalization, the electrostatic potential \(u\) associated with a prescribed charge density \(\rho\) satisfies the Born--Infeld equation
\begin{equation}\label{tradition}
	-\operatorname{div}\left(\frac{Du}{\sqrt{1 - |Du|^{2}}}\right) = \rho. 
\end{equation}

Born--Infeld-inspired modifications of gravity have also been investigated in stellar structure, compact objects, cosmological singularities, black holes, and wormholes \cite{BELTRAnJIMEnEZ20181,Casanellas_2012,PhysRevLett.111.071301,PhysRevD.93.104054}.

From the mathematical viewpoint,
the equation is related to spacelike graphs with prescribed Lorentzian mean curvature in Lorentz--Minkowski space \(\mathbb{L}^{N+1} = \mathbb{R}\times\mathbb{R}^N\), equipped with the metric \(ds^2 = -dx_0^2 + \sum_{j=1}^n dx_j^2\) \cite{WOS:A1980KV50800001}.

\begin{definition}
For $u \in W^{1,\infty}(\Omega)$, its graph is called
	
	\begin{itemize}
\item weakly spacelike if $|Du| \leq 1$ holds on $\Omega$;
\item spacelike if $|u(x) - u(y)| < |x - y|$ for any $x, y \in \Omega, x \neq y$ such that the closed segment $\overline{xy}$ joining \(x\) and \(y\) is contained in $\Omega$;
\item strictly spacelike if $u \in C^1(\Omega)$ and $|Du| < 1$ on $\Omega$.
	\end{itemize}
\end{definition}
\begin{definition}
Let $u:\Omega\to\mathbb{R}$ satisfy $|Du|\le 1$ almost everywhere. If $x,y\in\overline{\Omega}$, $x\neq y$, and $|u(y)-u(x)|=|x-y|$, then $\overline{xy}$ is called a light segment of $u$.
\end{definition}

On a bounded domain,

a spacelike graph of Lorentzian mean curvature \(\rho\) with boundary values \(\phi\) is described by a solution \(u : \overline{\Omega} \to \mathbb{R}\) of the Dirichlet problem
\begin{equation}\label{eq:bid}
	\begin{cases} 
		-\operatorname{div} \left( \dfrac{Du}{\sqrt{1 - |Du|^2}} \right) = \rho &  \Omega \subset \mathbb{R}^N , \\[1em]
		u = \phi & \partial\Omega ,
	\end{cases}
\end{equation}

In 1979, Flaherty \cite{5e5a24ac-4f55-3738-acec-10690ed03d3f} considered $\rho=0$ on a domain satisfying $\partial\Omega\in C^2$ and a nonnegative mean curvature condition, with $\phi\in C^2(\overline{\Omega})$ and $\max_{\overline{\Omega}}|D\varphi|<1$. He obtained a solution $u\in C^2(\overline{\Omega})$ satisfying $\max_{\overline{\Omega}}|Du|\le\gamma<1$ throughout $\overline{\Omega}$, where $\gamma$ depends only on $\phi$ and $\Omega$.

In 1982, Bartnik and Simon \cite{RN58, RN60} approached \eqref{eq:bid} with $\rho \neq 0$
by minimizing the action
\begin{equation}\label{eq:bii}
I_\rho(v)=\int_\Omega \left(1 - \sqrt{1 - |Dv|^2}\right) \mathrm{d}x - \langle \rho, v \rangle 
\end{equation}
in a variational formulation. For prescribed boundary values $\phi\in C(\partial\Omega)$, define
\begin{equation}
\mathcal{Y}_\phi(\Omega)=\{u\in W^{1,\infty}(\Omega):|Du|\le 1\ \text{a.e.},\ u|_{\partial\Omega}=\phi\}. 
\end{equation}
If the boundary values admit a spacelike extension $\bar\phi\in\mathcal{Y}_\phi(\Omega)$ on $\Omega$, then, for every $\rho\in L^\infty(\Omega)$, equation \eqref{eq:bid} admits a strictly spacelike weak solution $u\in C^1(\Omega)\cap W^{2,2}_{\rm loc}(\Omega)$.
\begin{definition}
Let $\Omega \subset \mathbb{R}^n$ be a bounded domain and $\rho \in W^{1,\infty}(\Omega)^*$. We call $u \in \mathcal{Y}_\phi(\Omega)$ a weak solution of \eqref{eq:bid} if
	\begin{align*}
	&(i) \displaystyle w \coloneq \frac{1}{\sqrt{1 - |Du|^2}} \in L_{\text{loc}}^1(\Omega);\\
	&(ii)\displaystyle \int_{\Omega} \frac{Du \cdot D\eta}{\sqrt{1 - |Du|^2}} \,dx = \langle \rho,  \eta \rangle, \qquad \forall \eta \in \operatorname{Lip}_c(\Omega).
	\end{align*}
Given a subdomain $\Omega' \subset \Omega$, if $w \in L_{\text{loc}}^1(\Omega')$ and (ii) holds for every $\eta \in \operatorname{Lip}_c(\Omega')$, then $u$ is a weak solution of \eqref{eq:bid} in $\Omega'$.
\end{definition}

When \(\rho\) is a finite sum of point charges, namely \(\rho = \sum_{j=1}^k a_j \delta_{x_j}, x_j \in \Omega\),
and $(\Omega, d_\Omega)$ has a compact completion while $\phi \in C(\partial\Omega)$ admits a spacelike extension,
Klyachin and Miklyukov \cite{klyachin1995existence}
proved the following statement.
For a fixed $k$-tuple $\mathscr{P} = (x_1, \ldots, x_k) \subset \Omega$,
there exists a constant $M_N(\phi,\mathscr{P})$
such that, for every $a \coloneq (a_1, \ldots, a_k) \in \mathbb{R}^k$ satisfying $|a| < M_n(\phi,\mathscr{P})$, the corresponding minimizer $u_\rho$ is strictly spacelike in $\Omega \setminus \{x_j\}$ and solves \eqref{eq:bid}.
In particular, if $N=2$, then $M_2(\phi,\mathscr{P}) = +\infty$, so there is no smallness restriction on the charges. If $N\ge 3$, a positive lower bound for this constant is determined by the zero-mean-curvature solution, the location of $\{x_1, \ldots, x_k\}$, and the geometry of the domain.

Byeon et al. \cite{WOS:001150613900001} subsequently extended the study of \eqref{eq:bid} from bounded sources $L^\infty$ and finitely many point charges to general Radon measure sources. They also related the weak solvability of the minimizer \(u_\rho\) of \eqref{eq:bii} to the presence of light segments.

For every \( N \geq 3 \) and \( \ell \in \{1, \ldots, N-2\} \),
one can choose $\rho=\rho_u\in L^q(\mathbb{R}^N)$ satisfying $q<N-\ell$
and a nontrivial function \( u \in C_c^\infty(\mathbb{R}^N) \) that weakly solves \eqref{eq:bid} with
\(\phi \equiv 0\).

\begin{definition}
If \(\Omega\subset\mathbb{R}^N\) is open, let \(\mathcal{M}(\Omega)\) denote the space of finite Borel measures on \(\Omega\), equipped with the total variation norm \(\|\cdot\|_{\mathcal{M}(\Omega)}\).
\end{definition}
Let \( \Sigma \Subset \Omega \) be compact and satisfy \( \mathcal{H}_\delta^1(\Sigma) = 0 \), where \(\mathcal{H}_{\delta}^{k}\) is the \(k\)-dimensional Hausdorff measure associated with the Euclidean distance \(\mathrm{d}_\delta\).
For boundary values $\phi\in C(\partial\Omega)$ admitting a spacelike extension, suppose that \( \rho \in \mathcal{M}(\Omega) \) decomposes as
\[
\rho = \rho_S + \rho_{AC}, \quad \quad \text{where}
\begin{cases} 
	\operatorname{supp} \rho_S \subset \Sigma,  \\ 
	\rho_{AC} \in L^1(\Omega) \cap L_{\text{loc}}^2(\Omega \setminus \Sigma). 
\end{cases}
\]
If $N=2$, the minimizer $u_\rho\in\mathcal{Y}_\phi(\Omega)$ weakly solves the equation in $\Omega$ and its graph contains no light segments.
If, in addition, $\Omega'\Subset\Omega\setminus\Sigma$ and $\rho\in L^\infty(\Omega')$,
then, for some $\alpha>0$, the minimizer $u_\rho$ is strictly spacelike in $\Omega'$ and belongs to $C^{1,\alpha}(\Omega')$.
When \( N \geq 3 \), let \( K_\phi^\rho \) denote the light-segment set of the minimizer \( u_\rho \in \mathcal{Y}_\phi(\Omega) \). Then $u_\rho$ weakly solves \eqref{eq:bid} in $\Omega\setminus K_\phi^\rho$.
Moreover, if \( K_\phi^\rho \cap (\partial\Omega \cup \Sigma) = \emptyset \), then \( u_\rho \) weakly solves \eqref{eq:bid} throughout \( \Omega \).

Concerning the relation between source integrability and light segments, De Filippis and Mingione \cite{DeFilippisMingione2026} used nonlinear potential theory to study the Born--Infeld equation with rough sources.
They extended the Bartnik--Simon anti-peeling result under suitable local Lorentz integrability assumptions. In the whole space, they combined dual-space assumptions with local integrability to obtain sufficient conditions for the absence of light segments and weak solvability of minimizers.
They also established quantitative strict spacelikeness. More precisely, suppose that \(N\ge3\),
and let \(\mathcal D^{1,2}(\mathbb R^N)\) be the completion of \(C_c^\infty(\mathbb R^N)\) with respect to the \(L^2\) norm of the gradient. If
\begin{equation}\label{eq:intro-DFM-source}
	\rho\in\mathcal D^{1,2}(\mathbb R^N)^*\cap L(N,1)(\mathbb R^N),
\end{equation}
and \(u_\rho\) minimizes the whole-space action over the admissible class
\(\{u\in\mathcal D^{1,2}(\mathbb R^N):\|Du\|_\infty\le1\}\), then \(u_\rho\) weakly solves \eqref{tradition} and satisfies
\begin{equation}\label{eq:intro-DFM-spacelike}
	\begin{aligned}
		\left\|\frac{1}{\sqrt{1-|Du_\rho|^2}}\right\|_{L^\infty(\mathbb R^N)}
		&\le \exp\{C_N[\rho]_{N,1}\},\\
		\|Du_\rho\|_{L^\infty(\mathbb R^N)}^2
		&\le 1-\exp\{-2C_N[\rho]_{N,1}\}<1.
	\end{aligned}
\end{equation}
Here \(L(N,1)\) denotes a Lorentz space, \([\rho]_{N,1}\) is its norm, and \(C_N>0\) depends only on the dimension. Moreover, \(Du_\rho\) is continuous on \(\mathbb R^N\).
Thus the corresponding assumptions identify minimizers as weak solutions and quantitatively separate the gradient from the singular boundary \(|Du|=1\).

If \(\Omega = \mathbb{R}^N\) and \( N \geq 3 \),
taking \(\rho = f(u)\) leads to the equation
\begin{equation}\label{eq:BI}
	\operatorname{div}\left(\frac{Du}{\sqrt{1 - |Du|^{2}}}\right) + f(u) = 0. 
\end{equation}
The results of Cheng--Yau and Ecker \cite{7f30351a-d61c-364e-871f-9ea19b79b56d,RN62}
on maximal graphs, including possible hyperplane profiles $u$ of slope $1$ when \(\rho = 0\), motivate the choice of a suitable function space for \(u\); the boundary condition in \eqref{eq:bid} is then replaced by \(u(x) \to 0\) as \(|x| \to \infty\).

Using a shooting method, Azzollini \cite{AZZOLLINI20142086} proved that, if \(f(u)\) satisfies
\[
\begin{cases}
	f(0)=0,\\[2pt]
f \text{ is locally Lipschitz on } [0,\infty), \\[2pt]
\lim\limits_{s\to\alpha^+}\dfrac{f(s)}{s-\alpha}>0,\quad\text{where } \alpha \coloneq \inf\{\xi>0 : f(\xi)\ge 0\}>0,\\[2pt]
\exists\gamma>0 \text{ such that } F(\gamma)\coloneq\int_0^\gamma f(s)\,ds>0 \\[2pt]
f(\xi)>0 \text{ on } (\alpha,\xi_0],\quad\xi_0 \coloneq \inf\{\xi>0 : F(\xi)>0\},
\end{cases}
\]
then equation \eqref{eq:BI} admits a radially decreasing ground state.

To work beyond radial symmetry and accommodate more general sources, Bonheure et al. \cite{RN60} introduced the Banach space
\begin{equation}\label{y}
	\mathcal{Y}(\mathbb{R}^N) \coloneq \overline{C_c^\infty(\mathbb{R}^N)}^{\|\cdot\|_{\mathcal{Y}}}, \quad \|v\|_{\mathcal{Y}} \coloneq \max \left\{ \|Dv\|_2, \|Dv\|_{p_1} \right\}.
\end{equation}
where \(p_1\) is fixed, $p_1\in(n,\infty)$, and, when $N=3$, one additionally assumes $p_1\ge 2^*=\frac{2N}{N-2}$.
In this space, if \(\rho \in L_{\mathrm{loc}}^{\infty}(\mathbb{R}^N)\), the minimizer \(u_\rho\) is locally strictly spacelike and \(u_\rho \in C_{\mathrm{loc}}^{1,\alpha}(\mathbb{R}^N)\) for some \(\alpha>0\).
For $\rho \in C_0^\alpha(\mathbb{R}^N) \cap L^1(\mathbb{R}^N)$,
\begin{equation*}
	\tau^{(1)}(x) \coloneq -\frac{N}{|S^{N-1}|} \nabla \int_{\mathbb{R}^N} \frac{\rho(y)}{|x-y|^{N-2}} \, \mathrm{d}^N y,
\end{equation*}
Carley and Kiessling \cite{RN66}, under the smallness condition $\| \tau^{(1)} \|_{1,\alpha} < (2^{2/3} - 1)^{3/2}$,
gave a constructive proof of the existence of a classical solution.

Byeon et al. \cite{WOS:001150613900001} then considered the closed convex subset of \eqref{y} given by
\begin{equation}
	\mathcal{Y}_0(\mathbb{R}^N) \coloneq \left\{ v \in \mathcal{Y}(\mathbb{R}^n) : \|Dv\|_\infty \leq 1 \right\},
\end{equation}
On this set, \(I_{\rho}\) admits a unique minimizer \(u_{\rho} \in \mathcal{Y}_0(\mathbb{R}^N)\).
If $\rho \in (L^1(\mathbb{R}^N) + L^p(\mathbb{R}^N)) \cap L_{\text{loc}}^2(\mathbb{R}^N)$ for some $p \in (1, 2_*](\text{where }2_*=\frac{2N}{N+2})$, then the minimizer $u_\rho$ weakly solves \eqref{eq:BI}.
If, moreover, $\rho \in L^\infty(\Omega')$, then $u_\rho$ is strictly spacelike in $\Omega' \Subset \Omega \setminus (\Sigma \cup K^\rho)$ and $u_\rho \in C_{\text{loc}}^{1,\alpha}(\Omega')$ for some $\alpha > 0$.

For the nonlinear source \(f(u)\), Byeon et al. \cite{Byeon2026} introduced a mass parameter $m \ge 0$. Suppose that $f \in C(\mathbb{R})$ and that there exists $ t_0 > 0$ such that
\begin{equation}
	F(t_0) \equiv \int_0^{t_0} f(t)\,\mathrm{d}t > 0, 
\end{equation}
Assume also that, when $m = 0$,
\begin{equation}
	-\infty < \liminf_{s\to 0}\frac{f(s)}{|s|^{2^*-2}s} \le \limsup_{s\to 0}\frac{f(s)}{|s|^{2^*-2}s} = 0,
\end{equation}
and, when $m > 0$,
\begin{equation}
	-\infty < \liminf_{s\to 0}\frac{f(s)}{|s|^{p-2}s} \le \limsup_{s\to 0}\frac{f(s)}{|s|^{p-2}s} = -m, 
\end{equation}
for some $p\in[2,2^*]$.
Define
\begin{equation}
	\mathcal{D}_{m,p}^{1,2}(\mathbb{R}^N) \coloneq \overline{C_c^\infty(\mathbb{R}^N)}^{\|\cdot\|_{\mathcal{D}_{m,p}^{1,2}}}, \qquad
	\|\phi\|_{\mathcal{D}_{m,p}^{1,2}} \coloneq \sqrt{\|D\phi\|_2^2 + m\|\phi\|_p^2}.
\end{equation}
Then \eqref{eq:BI} admits a positive radial weak solution $u \in \mathcal{D}_{m,p}^{1,2}(\mathbb{R}^N)$.
If $f$ is odd, \eqref{eq:BI} admits infinitely many distinct radial sign-changing weak solutions $\{u_k\}_k \subset \mathcal{D}_{m,p}^{1,2}(\mathbb{R}^N)$.
If, in addition, $N,d$ satisfy
\[
N \geq 4,\qquad 2 \leq d \leq \frac{N}{2},\qquad \text{and}\quad N - 2d \neq 1,
\]
there are also infinitely many distinct nonradial sign-changing weak solutions $\{u_k\}_k \subset \mathcal{D}_{m,p}^{1,2}(\mathbb{R}^N)$.
Furthermore, every weak solution $v \in \mathcal{D}_{m,p}^{1,2}(\mathbb{R}^N)$ of \eqref{eq:BI} satisfies $v \in W_{\mathrm{loc}}^{2,q}(\mathbb{R}^N)$ for all $q \in [2,\infty)$, as well as $\|Dv\|_\infty < 1$ in $\mathbb{R}^N$.

For ground states, Bieganowski et al. \cite{BieganowskiIkomaMederski2026} studied the autonomous equation \eqref{eq:BI} in both the zero-mass and positive-mass cases.
For example, in the normalized positive-mass case, write \(f(s)=-s+g(s)\), where \(g\in C^1(\mathbb R)\) satisfies
\[
\begin{cases}
	g(0)=g'(0)=0,\\[2pt]
s\longmapsto\dfrac{g(s)}{|s|}\ \text{is strictly increasing on }\mathbb R\text{, with its continuous extension at }s=0,\\[4pt]
	\displaystyle\lim_{|s|\to\infty}\dfrac{g(s)}{s}=+\infty.
\end{cases}
\]
Fix an integer \(j_0\ge1\) such that \(2j_0>N\), and define the admissible class and the energy by
\begin{equation}\label{eq:intro-BIM-energy}
	\begin{aligned}
		\mathcal X_{\mathrm a}
		&=\{u\in H^1(\mathbb R^N):Du\in L^{2j_0}(\mathbb R^N),\ \|Du\|_\infty\le1\},\\
		\mathcal J(u)
		&=\int_{\mathbb R^N}\left(1-\sqrt{1-|Du|^2}\right)\,dx
		-\int_{\mathbb R^N}F(u)\,dx,
	\end{aligned}
\end{equation}
where \(F(s)=\int_0^s f(t)\,dt\). Analyzing the gradient-preserving scaling \(u_t(x)=t u(x/t)\) leads to the Pohozaev constraint
\begin{equation}\label{eq:intro-BIM-Pohozaev}
	\begin{aligned}
		\mathcal P(u)&=\mathcal J(u)-\frac1N\int_{\mathbb R^N}f(u)u\,dx,\\
		\mathcal M_{\mathrm a}&=\{u\in\mathcal X_{\mathrm a}\setminus\{0\}:\mathcal P(u)=0\},
	\end{aligned}
\end{equation}
The associated minimization problem is attained, and its ground-state level satisfies
\begin{equation}\label{eq:intro-BIM-ground-level}
	0<c_{\mathrm{gs}}:=\inf_{u\in\mathcal M_{\mathrm a}}\mathcal J(u)=c_{\mathrm{MP}},
\end{equation}
where \(c_{\mathrm{MP}}\) is the corresponding mountain-pass level. A profile decomposition adapted to the nonsmooth functional yields strong convergence of minimizing sequences after suitable translations, and the limit is a classical ground state of \eqref{eq:BI}.
This construction does not restrict the minimization problem to radial functions in advance. Constant sign and the existence of a radial ground state are established afterward.
Under the corresponding oddness and other assumptions on the nonlinearity, the authors also constructed nonradial solutions for \(N\ge4\).
These results develop the variational theory of autonomous Born--Infeld equations from the least-energy viewpoint.

\subsection{Nonsmooth critical point theory}

Mathematical models such as the Laplace, wave, and nonlinear Schr{\"o}dinger equations
describe equilibria, evolution, or standing waves in physical systems.
Nonlinearity, singularities, and unbounded domains, however,
often prevent the explicit construction of solutions.

Variational principles associate an energy or action functional with the equation, and equilibria correspond to critical points, which may be minima or saddle points.
They therefore convert the differential equation into a critical point problem \cite{Basdevant2007}.
Under appropriate regularity and variational correspondence assumptions, these critical points yield weak or generalized solutions.

In the smooth setting, the mountain pass theorem is a basic tool for constructing nontrivial solutions.
Let \(X\) be a Banach space and suppose that \(I\in C^1(X,\mathbb{R})\) satisfies
\begin{enumerate}[label=(\roman*)]
	\item \(I(0)=0\);
\item there exist \(\rho>0\) and \(\alpha>0\) such that
	\(\inf_{\|u\|=\rho}I(u)\ge\alpha;\)
\item there exists \(e\notin B_\rho(0)\) such that
	\(I(e)\le0.\)
\end{enumerate}
Then \(I\) is said to have mountain pass geometry \cite{Willem1996}.
Geometrically, the sphere \(\partial B_\rho\) separates the origin \(0\) from \(e\) by an energy barrier that any path from \(0\) to \(e\) must cross.
Define the class of continuous paths joining \(0\) to \(e\) by
\[
\Gamma
=
\bigl\{
\gamma\in C([0,1],X):
\gamma(0)=0,\ \gamma(1)=e
\bigr\},
\]
Minimizing the largest energy along these paths leads to the minimax value
\[
c
=
\inf_{\gamma\in\Gamma}
\max_{t\in[0,1]}I(\gamma(t)).
\]
To pass from this geometric construction to a critical point, one introduces a Palais--Smale sequence at level \(c\), abbreviated as a \((PS)_c\) sequence.
Thus \(\{u_n\}\subset X\) satisfies
\(I(u_n)\to c,\)
and \(I'(u_n)\to0\) in \(X^*\)
\cite{PALAIS1966115,34c52b66-b753-302a-96ef-9c37edd724b2}. 

Under mountain pass geometry, such a sequence can be obtained by a deformation lemma \cite{https://doi.org/10.1002/cpa.3160330507} or Ekeland's variational principle \cite{EKELAND1974324}.
If the sequence diverges or oscillates without a convergent subsequence, it need not yield a critical point.
An additional compactness requirement, the Palais--Smale condition or PS condition, is therefore imposed.
It requires every \((PS)_c\) sequence to have a strongly convergent subsequence and is central to deformation arguments and the existence \cite{EKELAND1974324} and multiplicity \cite{rabinowitz1986minimaxBI253} of critical points.
For many problems on bounded domains, compact Sobolev embeddings yield this condition, which has consequently found extensive applications \cite{Willem1996}.
The mountain pass theorem then gives \(u\in X\) such that
\[
I(u)=c,
\qquad
I'(u)=0.
\]
Hence \(c\) is a critical value of \(I\); see Jabri \cite{Jabri2011}.
This completes the passage from minimax geometry to a critical point.

For the Born--Infeld model, the variational functional associated with \eqref{tradition} is
\[
I_\rho(v) 
\coloneq 
\int_\Omega 
\left(1 - \sqrt{1 - |Dv|^2}\right) \mathrm{d}x 
- 
\langle \rho, v \rangle .
\]
The derivative of the integrand \(1 - \sqrt{1 - |Dv|^2}\) becomes unbounded as \(|Dv| = 1\).
At points where \(|Dv| = 1\), the G{\^a}teaux derivative may fail to exist as a finite quantity.
Consequently, the classical \(C^1\) mountain pass theorem does not apply directly to the Born--Infeld functional.

For an energy that is not differentiable everywhere, critical points can instead be defined by a variational inequality, allowing the study of existence, multiplicity, and qualitative properties in a nonsmooth setting.
Such methods have applications to minimal hypersurfaces, multiplicity for semilinear and quasilinear elliptic problems, finite elasticity, periodic Hamiltonian trajectories, and solitary waves of semilinear wave equations \cite{ref1}.

Szulkin\cite{SZULKIN198677}
considered functionals split into a convex part and a \(C^1\) part. Let \(X\) be a Banach space and write
\[
I(u)=\Psi(u)-\Phi(u),
\]
where \(\Psi\) and \(\Phi\) satisfy
\begin{itemize}
\item[\((\Psi_1)\)] \(\Psi : X \to [0, \infty]\) is lower semicontinuous and convex, and \(\Psi(0) = 0\);
	
	\item[\((\Phi_1)\)] \(\Phi \in C^1(X, \mathbb{R})\).
\end{itemize}
\begin{definition}
A point \(u\) is critical for \(I\) if \(u\in D(\Psi)=\{u:\Psi(u)<+\infty\}\) and
	\[
	\Psi(v)-\Psi(u)-\Phi'(u)(v-u)\ge0,
	\qquad
	\forall v\in X.
	\]
\end{definition}
\begin{definition}
A sequence \(\{u_j\}\subset D(\Psi)\) is a PS sequence for \(I\) at level \(c\) if \(I(u_j)\to c\) and there exists \(\varepsilon_j\to0^+\) such that
	\[
	\Psi(v)-\Psi(u_j)-\Phi'(u_j)(v-u_j)
	\ge 
	-\varepsilon_j\|v-u_j\|,
	\qquad
	\forall v\in X.
	\]
If every such sequence has a strongly convergent subsequence, then \(I\) satisfies the PS condition.
\end{definition}
Under this nonsmooth PS condition, the relevant sequences possess strongly convergent subsequences,
which permit energy descent in deformation arguments.
This replaces the direct use of a classical gradient flow for a nondifferentiable functional
and underlies nonsmooth extensions of mountain pass and saddle point theorems.

For second- and higher-order elliptic operators \(L\), Chang \cite{CHANG1981102}
introduced the Clarke generalized gradient
to treat functionals with weaker regularity. Consider
\[
Lu(x) \in 
\bigl[ 
\underline{\phi}\bigl(x,u(x)\bigr),\; 
\overline{\phi}\bigl(x,u(x)\bigr) 
\bigr], 
\qquad 
x\in\Omega,
\qquad 
u|_{\partial\Omega}=0,
\]
where \(\underline{\phi},\overline{\phi}\) are the essential lower and upper envelopes of the nonlinearity \(\phi(x,t)\) near \(t\), namely
\[
\underline{\phi}(x,t)
=
\lim_{\delta\to0}
\underset{|s-t|<\delta}{\operatorname{ess\,inf}}\;
\phi(x,s),
\qquad
\overline{\phi}(x,t)
=
\lim_{\delta\to0}
\underset{|s-t|<\delta}{\operatorname{ess\,sup}}\;
\phi(x,s).
\]
Under suitable local boundedness and sublinear or superlinear growth assumptions on \(\phi(x,t)\),
critical points can be obtained in the sense of the Clarke generalized gradient,
yielding nontrivial solutions to the corresponding differential inclusion.

On the whole space \(\mathbb{R}^{N}\), translation invariance
may allow PS sequences to converge only weakly or to lose mass at infinity. At critical exponents or under geometric invariances,
dilations provide a further source of noncompactness.
In 1988, Struwe \cite{10.1007/BF02392272} introduced a parameter and exploited monotonicity to address noncompactness in a conformally invariant geometric problem.
In 1999, Jeanjean \cite{Jeanjean_1999} developed an abstract version of this argument.
More precisely, let \(X\) be a Banach space,
\[
I_\lambda(u) = \Psi(u) - \lambda \Phi(u), 
\qquad 
\lambda \in J \subset \mathbb{R}^+, 
\qquad 
\Phi(u) \ge 0,
\]
and suppose that \(\Psi(u) \to +\infty\) or \(\Phi(u) \to +\infty\) as \(\|u\| \to \infty\).
If each \(I_\lambda\) has mountain pass level \(c_\lambda> \max\{I_\lambda(v_1), I_\lambda(v_2)\}\), as in Definition \ref{def:mountain} below,
then for almost every \(\lambda \in J\) there is a bounded PS sequence \(\{v_n\} \subset X\) such that
\[
I_\lambda(v_n) \to c_\lambda, 
\qquad 
I_\lambda'(v_n) \to 0 
\quad 
\text{in } X^*.
\]
This replaces the difficult direct boundedness problem by critical or approximate critical points of nearby functionals \(I_{\lambda_n}\).
When \(\lambda_n\to1\) and these points are bounded in \(X\), suitable bounds for \(\Phi,\Phi'\) on bounded sets show that they also form a PS sequence for the original functional \(I=I_1\).
In an application, Jeanjean worked in \(X=H^{1}(\mathbb{R}^{N})\) and considered
\[
-\Delta u+Ku=f(x,u),
\qquad 
x\in\mathbb{R}^{N},
\qquad 
K>0,
\]
with energy
\[
I(u)
=
\frac12\int_{\mathbb{R}^{N}}
\bigl(|\nabla u|^{2}+Ku^{2}\bigr)\,dx
-
\int_{\mathbb{R}^{N}}F(x,u)\,dx.
\]
The monotonicity trick and assumptions on the nonlinearity first produce a bounded PS sequence at the mountain pass level.
The concentration-compactness method of Lions \cite{Lions1984} then excludes vanishing, while a weak maximum principle gives nonnegativity and nontriviality of the limit.
Thus \(I\) has a nontrivial critical point, yielding a positive solution of the original equation.

In 2012, Squassina \cite{Squassina_2012} extended the argument to continuous functionals on Banach spaces, including geometric and nonlinear elliptic problems without \(C^1\) regularity.
In that framework, criticality is characterized
by vanishing weak slope rather than vanishing Fr{\'e}chet derivative, and PS sequences have convergent energy and weak slope tending to zero.
Squassina also combined the monotonicity trick with symmetrization.
If a Banach space \(X\) embeds continuously into \(V\) and carries an abstract symmetrization map,
one constructs bounded, almost symmetric PS sequences
and exploits compact embeddings in symmetric function spaces to recover compactness.

To address constrained nonsmooth problems and the divergence of higher minimax levels,
Byeon et al. \cite{Byeon2026} developed the monotonicity method for lower semicontinuous functionals of Szulkin type.
Unlike the continuous-functional framework, this setting allows the value \(+\infty\)
and naturally accommodates the Born--Infeld constraint \(|Du|\le1\).
In place of the global PS condition, one assumes boundedness of the following critical-point sets:
for every \(0 < a \leq b < \infty\) and \(c \in \mathbb{R}\),
\[
\mathcal{K}_{[a,b]}^c 
\coloneq 
\bigl\{ 
u \in X 
\mid 
\text{there exists } \lambda \in [a,b] ,
u \text{ is critical for } I_\lambda \text{ and } I_\lambda(u) \leq c
\bigr\}
\]
is bounded in \(X\).
This requires uniform bounds for energy-controlled families of genuine critical points over compact parameter intervals, rather than compactness of arbitrary PS sequences.
For many elliptic equations, such bounds follow from a Pohozaev-type identity.

Byeon et al. first establish uniform mountain pass geometry and define a minimax level for each parameter \(\lambda\).
Because \(\Psi\ge 0\), the minimax levels are monotone in the parameter, and the monotonicity trick produces bounded PS sequences for almost every parameter.
Under the relevant lower semicontinuity, weak stability, and bounded-PS compactness assumptions,
their limits are genuine critical points of the corresponding \(I_\lambda\).
One then lets the parameter tend to the target value \(\lambda=1\)
and uses boundedness of \(\mathcal{K}_{[a,b]}^c\) to obtain a critical point of the original functional \(I=I_1\).
For even functionals, symmetric minimax classes yield a sequence of levels and infinitely many critical points.
Applications to autonomous Born--Infeld equations give radial and nonradial solutions.

\subsection{The non-autonomous problem and main results}

We study the existence of solutions to the non-autonomous Born--Infeld equation
\begin{equation}\label{biff}
	\operatorname{div}\left( \frac{Du}{\sqrt{1 - |Du|^2}} \right) - V(x)u + |u|^{p-2}u = 0, \quad x \in \mathbb{R}^N, 
\end{equation}
where \(N\ge 3,\ 2<p<2^*\) and the potential \(V \in C^1(\mathbb{R}^N)\) satisfies
\[
V_0 \coloneq \inf_{x \in \mathbb{R}^N} V(x) > 0
\]
and
\begin{equation}\label{Qx}
	\lim_{|x| \to \infty} V(x) = +\infty.
\end{equation}
To construct negative-energy endpoints and obtain uniform bounds for critical points, we further assume
\begin{equation}\label{V-derivative-growth}
	|x \cdot DV(x)| \leq C_D V(x), \qquad C_D < p - 2.
\end{equation}

The associated functional \(I = \Psi(u) - \Phi(u)\) is defined by
\begin{equation}
	\Psi(u)=\begin{cases}
		\displaystyle
		\int_{\mathbb{R}^N} \bigl( 1 - \sqrt{1 - |Du|^2} \bigr) \, dx + \frac{1}{2} \int_{\mathbb{R}^N} V(x)u^2 \, dx,
		 & \|Du\|_\infty\le 1,\\[6pt]
		\infty, & \text{else},
	\end{cases}
\end{equation}
\begin{equation}
	\Phi(u)=\frac{1}{p} \int_{\mathbb{R}^N} |u|^p \, dx.
\end{equation}

We investigate the problem in the natural space \(Y\coloneq H^1_V(\mathbb R^N)\), where
\begin{equation}
	H_V^1(\mathbb{R}^N) = \left\{ u \in H^1(\mathbb{R}^N) : \int_{\mathbb{R}^N} V(x)u^2 \, dx < +\infty \right\},
\end{equation}
with norm
\begin{equation}
	\|u\|_{Y} = \left( \int_{\mathbb{R}^N} \bigl( |Du|^2 + V(x)u^2 \bigr) \, dx \right)^{1/2}.
\end{equation}
The admissible set is
\begin{equation}\label{tr-Dpsi}
	D_Y(\Psi)
	=
	\left\{
	u\in Y:\|Du\|_\infty\le 1
	\right\}.
\end{equation}
\begin{remark}
For the variational construction, we introduce an auxiliary space \(X\). Its definition, norm, and relation to \(Y\) are established in the next section.
\end{remark}

We shall verify the following properties of \(\Psi\) and \(\Phi\):

\begin{itemize}
    \item[$(\Phi_1)$] $\Phi \in C^1(X, \mathbb{R})$; 

\item[$(\Psi_1)$] $\Psi : X \to [0, \infty]$ is lower semicontinuous and convex, and $\Psi(0) = 0$;
	
\item[$(\Psi_2)$] $\Psi(u) \to \infty$ as $\|u\| \to \infty$;
	
\item[$(\Psi_3)$] if $\{u_j\}_{j=1}^\infty \subset D(\Psi)$ and $u \in D(\Psi)$ satisfy $u_j \rightharpoonup u$ in $X$ and $\Psi(u_j) \to \Psi(u)$, then $\|u_j - u\| \to 0$.
	
\end{itemize}
Write
\begin{equation}
	D(\Psi) \coloneq \{ u \in X \mid \Psi(u) < \infty \}, 
\end{equation}

To construct bounded PS sequences, we introduce the parameter-dependent family
\begin{definition}\label{def:critical}
Let $X$ be a Banach space and $I_\lambda = \lambda \Psi - \Phi$, where $\Psi,\Phi$ satisfies $(\Psi_1),(\Phi_1)$ and $\lambda \in \mathbb{R}^+$.
	\begin{itemize}
\item A point $u \in X$ is critical for $I_\lambda$ if, for every \(v \in X\),
		\begin{equation}\label{nonsm-critical}
			\lambda (\Psi(v) - \Psi(u)) - \Phi'(u)(v - u) \ge 0.  
		\end{equation}
\item A sequence $\{u_j\}_{j=1}^\infty \subset X$ is a Palais--Smale sequence for $I_\lambda$ at level $c \in \mathbb{R}$ if
		\begin{itemize}
\item[(i)] $I_\lambda(u_j) \to c$ as $j \to \infty$;
\item[(ii)] there exists $\{\varepsilon_j\}_{j=1}^\infty \subset \mathbb{R}^+$ satisfying $\varepsilon_j \to 0$ as $j \to \infty$ such that, for every \(v \in X\),
			\begin{equation}\label{PSineq}
				\lambda (\Psi(v) - \Psi(u_j)) - \Phi'(u_j)(v - u_j) \ge -\varepsilon_j \|v - u_j\|.
			\end{equation}
		\end{itemize}
	\end{itemize}
\end{definition}

For \(\lambda\) close to \(1\), the family \(I_\lambda = \lambda\Psi - \Phi\) has uniform mountain pass geometry: there exist positive \(\rho_0, \alpha_0, \varepsilon_0\) and \(e \in X\) such that
\begin{equation*} \label{eq:MP}
	\text{(MP)} \quad \left\{
	\begin{array}{l}
		\inf\limits_{\|u\|=\rho_0} I_{1-\varepsilon_0}(u) \equiv \alpha_0 > 0,\\
		 \max \{I_{1+\varepsilon_0}(0), I_{1+\varepsilon_0}(e)\} \leq 0, \qquad
		  \|e\| > \rho_0. 
	\end{array}
	\right.
\end{equation*}

\begin{remark}
The endpoint \(e\) is constructed in Lemma \ref{MP-Q} and remains fixed in all subsequent path classes.
\end{remark}

\begin{definition}\label{def:mountain}
For \(\lambda \in [1-\varepsilon_0, 1+\varepsilon_0]\), define the mountain pass level
	\begin{equation} \label{eq:c0_lambda}
		c_0(\lambda) \coloneq \inf_{\gamma \in \Gamma_0} \sup_{t \in [0,1]} I_\lambda(\gamma(t)),
	\end{equation}
where
	\begin{equation} \label{eq:Gamma0}
		\Gamma_0 \coloneq \{ \gamma \in C([0,1], X) \mid \gamma(0)=0, \ \gamma(1)=e \}.
	\end{equation}
\end{definition}

The family \(I_\lambda=\lambda\Psi-\Phi\) is also even, that is,
\begin{equation}\label{condition-E}
	\Psi(-u)=\Psi(u),
	\qquad
	\Phi(-u)=\Phi(u),
	\qquad
	u\in X.
	\tag{E}
\end{equation}
We shall verify the following symmetric mountain pass geometry:

for every \(k\in\mathbb N\), write
\[
D^k(r)
=
\left\{
\zeta\in\mathbb R^k:|\zeta|\le r
\right\},
\qquad
D^k=D^k(1).
\]
there exist \(\rho_0,\alpha_0,\varepsilon_0>0,\) and, for each \(k\in\mathbb N,\),
a continuous odd map \(\pi_{0,k}\in
C(\partial D^k,X)\) such that
\begin{equation*} \label{eq:SMP}
	\text{(SMP)} \quad \left\{
	\begin{array}{l}
		\inf_{\|u\|_X=\rho_0}
		I_{1-\varepsilon_0}(u)
		\ge\alpha_0>0,\\
		\min_{\zeta\in\partial D^k}
		\|\pi_{0,k}(\zeta)\|_X>\rho_0,\\
		\sup_{\zeta\in\partial D^k}
		I_{1+\varepsilon_0}
		\bigl(\pi_{0,k}(\zeta)\bigr)<0.
	\end{array}
	\right.
\end{equation*}
\begin{definition}\label{def:symmetric-mountain}
For every
	\(\lambda\in[1-\varepsilon_0,1+\varepsilon_0]\)
and \(k\in\mathbb N\), define
	\begin{equation}\label{definition-ck}
		c_k(\lambda)
		=
		\inf_{\gamma\in\Gamma_k}
		\sup_{\zeta\in D^k}
		I_\lambda(\gamma(\zeta)),
	\end{equation}
where
	\begin{equation}\label{definition-Gammak}
		\Gamma_k
		\coloneq
		\left\{
		\gamma\in C(D^k,X):
\gamma\text{ is odd and }
		\gamma|_{\partial D^k}=\pi_{0,k}
		\right\}.
	\end{equation}
\end{definition}

Under the preceding geometric and functional hypotheses, bounded PS sequences are obtained as follows.
\begin{lemma}[Existence of bounded PS sequences \cite{Byeon2026}]\label{Psexist}
Suppose that \(\Psi\) and \(\Phi\) satisfy
\((\Psi_1)\), \((\Psi_2)\), and \((\Phi_1)\), with \(\lambda \in (1-\varepsilon_0, 1+\varepsilon_0)\).
	\begin{itemize}
\item If (MP) holds and
\(c_0'(\lambda)\) exists, then \(I_\lambda\)
has a bounded PS sequence in \(X\)
\(\{u_j^{(0)}\}_{j=1}^\infty \subset D(\Psi) \), abbreviated as \(\{u_j\}\), satisfying
		\(I_\lambda(u_j)\to c_0(\lambda).\)
\item If (SMP) and (E) hold and
\(c_k'(\lambda)\) exists, then \(I_\lambda\) has a bounded PS sequence in \(X\)
\(\{u_j^{(k)}\}_{j=1}^\infty\subset D(\Psi)\) satisfying
		\(I_\lambda(u_j^{(k)})\to c_k(\lambda).\)
	\end{itemize}
\end{lemma}
In addition, suppose that
\begin{itemize}
\item[$(\Phi_2)$] every bounded sequence \(\{u_l\}_l\) in \(X\) has a weakly convergent subsequence \(\{u_{l_j}\}_j\) and a point \(u\in X\) such that
	\begin{equation}
		\liminf_{j\to\infty} \Phi'(u_{l_j})(v-u_{l_j})
		\ge \Phi'(u)(v-u),\qquad \forall v\in X.
	\end{equation}
\item[(IB)] for every \(0<a\le b<\infty\) and \(c\in\mathbb R\), the set
	\begin{equation*}
		\mathcal K_{[a,b]}^c :=
		\left\{ u\in X : \exists\,\lambda\in[a,b],\ 
u\text{ is critical for }I_\lambda\text{ and }I_\lambda(u)\le c
		\right\}
\text{ is bounded in }X.
	\end{equation*}
\end{itemize}
Then the following multiplicity result applies.
\begin{lemma}[Existence of critical points \cite{Byeon2026}]\label{th1.7}
Assume $(\Psi_1)$, $(\Psi_2)$, $(\Psi_3)$, $(\Phi_1)$, $(\Phi_2)$, (IB), (E), and (SMP). Then $I = I_1$ has infinitely many critical points $\{u^{(k)}\}_{k=1}^\infty$ in $X$ satisfying $I(u^{(k)}) = c_k(1) \to \infty$ as $k \to \infty$.
\end{lemma}
For the symmetric minimax levels with \(k\ge1\),
we apply \ref{th1.7} directly after verifying \((\Phi_2)\), (IB), (E), and (SMP), obtaining infinitely many critical points of \(I_1\).
For the mountain pass level \(c_0(\lambda)\), we construct parameter-dependent critical points by the monotonicity trick and pass to \(\lambda_n\downarrow1\) to obtain a nontrivial critical point of \(I_1\).

\begin{theorem}
For each \(n\in\mathbb N\), the strong limit \(u_n\in D(\Psi)\) is a critical point of \(I_{\lambda_n}\) at \(c_0(\lambda_n)\), namely
	\[
	I_{\lambda_n}(u_n)=c_0(\lambda_n),
	\]
and
	\[
	\lambda_n\bigl(\Psi(v)-\Psi(u_n)\bigr)
	-
	\Phi'(u_n)(v-u_n)
	\ge0,
	\qquad
	\forall v\in D(\Psi).
	\]
\end{theorem}

Combining this construction with regularity yields the following results.

\begin{theorem}
The original problem \eqref{biff} admits a nontrivial positive weak solution \(w\in Y\).
\end{theorem}
Applying Lemma \ref{th1.7} also gives
\begin{theorem}
The original problem \eqref{biff} admits infinitely many distinct nontrivial weak solutions \(\{u^{(k)}\}_{k=1}^{\infty}
	\subset Y.\)
\end{theorem}
All the weak solutions obtained above enjoy the following regularity.
\begin{theorem}
The positive weak solution \(w\) of Theorem \ref{solution-positive}
and all weak solutions of Theorem \ref{solution-infty}
\(u^{(k)}\) satisfy, for every \(q>N\),
	\[
	w,u^{(k)}\in W_{\mathrm{loc}}^{2,q}(\mathbb R^N),
	\]
and, in \(\mathbb R^N\),
	\[
	|Dw|<1,
	\qquad
	|Du^{(k)}|<1.
	\]
Moreover, for some \(\alpha\in(0,1)\),
	\[
	w,u^{(k)}\in
	C_{\mathrm{loc}}^{2,\alpha}(\mathbb R^N).
	\]
Thus \(w\) is a positive classical solution of the original equation,
and \(\{u^{(k)}\}_{k=1}^{\infty}\) are infinitely many
nontrivial classical solutions.
\end{theorem}

\subsection{Organization and notation}
The paper is organized as follows.
Section 2 introduces the function spaces, properties of PS sequences, and auxiliary estimates.
Section 3 verifies the basic functional assumptions and evenness,
establishes uniform mountain pass and symmetric mountain pass geometry, and constructs bounded PS sequences at \(c_0(\lambda)\) by the monotonicity trick.
Compact embeddings and condition \((\Psi_3)\) yield strong convergence to critical points of the corresponding parameter-dependent functionals.
Section 4 identifies critical points for every \(\lambda>0\) with weak solutions and establishes local \(W^{2,q}\) regularity and strict spacelikeness.
Section 5 first derives an identity that provides uniform \(X\) bounds for critical points,
then proves the existence of a positive weak solution and infinitely many distinct weak solutions,
and finally applies local Schauder theory to obtain classical regularity for all these solutions.

For convenience, we introduce the following notation for the integrals in the functional.
Set \(G(\xi)=1-\sqrt{1-|\xi|^2}\),
so that
\begin{equation}
	\int_{\mathbb{R}^N} \bigl( 1 - \sqrt{1 - |Du|^2} \bigr) \, dx=\int_{\mathbb R^N}G(Du)\,dx
\end{equation}
Also write
\begin{equation}\label{ku}
	\mathcal K(u)
	=
	\int_{\mathbb R^N}G(Du)\,dx.
\end{equation}
For the regularity analysis, we express the parameter-dependent equation as
\begin{equation}\label{hceq}
	\operatorname{div}\left(\frac{Du}{\sqrt{1-|Du|^2}}\right)
	-
	V(x)u
	+
	\lambda^{-1}|u|^{p-2}u
	=
	0,
	\qquad
	x\in\mathbb R^N.
\end{equation}
Write
\begin{equation}\label{rholu}
	\rho_{\lambda,u}(x)
	=
	\lambda^{-1}|u|^{p-2}u
	-
	V(x)u,
\end{equation}
which can be rewritten as
\begin{equation}\label{hceqjx}
	\operatorname{div}\left(\frac{Du}{\sqrt{1-|Du|^2}}\right)
	+
	\rho_{\lambda,u}
	=
	0.
\end{equation}
If \(u\in C^2(\mathbb R^N)\) and \(|Du|<1\), then \(M = \{ (u(x), x) : x \in \mathbb{R}^N \}\)
is a strictly spacelike graph with Lorentzian mean curvature \(H_\lambda\) satisfying
\begin{equation}\label{H}
	H_\lambda
	=
	\operatorname{div}\left(\frac{Du}{\sqrt{1-|Du|^2}}\right)
	=
	V(x)u-\lambda^{-1}|u|^{p-2}u.
\end{equation}
This identity is justified in Section \ref{compute}.
When \(\lambda=1\), the mean curvature for the original problem is
\[
H
\coloneq H_1=
V(x)u-|u|^{p-2}u.
\]

\begin{remark}\label{Hrelax}
In the application, regularity gives, for every \(q>N\) in \(\mathbb R^N\) and \(u\in W_{\mathrm{loc}}^{2,q}(\mathbb R^N),\),
\(u\in C_{\mathrm{loc}}^{1,\alpha}(\mathbb R^N)\), where $\alpha$=1-$\frac{N}{q}$.
Together with the local integrability \(|Du|<1\) supplied by strict spacelikeness, this permits the mean-curvature identity in the distributional sense; see Theorem \ref{thm:critical-to-weak-regular}.
\end{remark}

\section{Functional setting and preliminary lemmas}

We first construct the auxiliary space \(X\) and establish its relation to the natural space \(Y\)
to prepare the variational argument in \(X\). We then record properties of PS sequences, the nonsmooth critical point formulation, and the scaling estimates needed below.

\subsection{Construction and properties of the spaces}
Fix \(p_1>\max\{N,p\}\).
On \(C_c^\infty(\mathbb R^N)\), define the norm
\[ 
\|u\|_X = \left( \|u\|_Y^2+\|Du\|_{p_1}^2 \right)^{1/2}. 
\] 
that is,
\begin{equation}\label{Xfanshu}
	\|u\|_X = \left( \int_{\mathbb R^N} \left( |Du|^2+V(x)u^2 \right)\,dx + \|Du\|_{p_1}^2 \right)^{1/2}. 
\end{equation}
Define
\begin{equation}\label{fzr} 
	X = \overline{C_c^\infty(\mathbb{R}^N)}^{\|\cdot\|_X}. 
\end{equation} 
Since \(V(x)\ge V_0>0\), the embedding \(Y\hookrightarrow H^1(\mathbb R^N)\) is continuous. Indeed, \[
\int_{\mathbb R^N}u^2\,dx \le \frac1{V_0} \int_{\mathbb R^N}V(x)u^2\,dx, 
\] 
and hence
\[
 \|u\|_{H^1}^2 \le C_V\|u\|_Y^2. 
 \] 
It is immediate that
 \[ 
 X\subset Y, \qquad \|u\|_Y\le \|u\|_X. 
 \]

We now relate this space to the intersection \(Y\cap W^{1,p_1}(\mathbb R^N)\). We first recall the following Gagliardo--Nirenberg inequality.

\begin{lemma}[Gagliardo--Nirenberg inequality \cite{LiuWangGNBestConstant}]\label{lem:Gn}
Let \(d \ge 1\), \(s > d\), and \(q \ge 0\), and define
	\[
	\theta = \frac{s d}{d s + (s - d)(q + 1)} \in (0,1).
	\]
There is an optimal constant \(C_{q,\infty,s} > 0\) such that every function
	\(u \in \mathcal{D}^{s,q+1}(\mathbb{R}^d)\)
satisfies
	\[
	\|u\|_{L^{\infty}} \le C_{q,\infty,s} \|u\|_{L^{q+1}}^{1-\theta} \|\nabla u\|_{L^s}^{\theta}.
	\]
where
	\[
	D^{s,q+1}(\mathbb{R}^d) = \overline{C_c^\infty(\mathbb{R}^d)}^{\|\cdot\|_{D^{s,q+1}}},
	\qquad
	\|u\|_{D^{s,q+1}} = \|u\|_{q+1} + \|Du\|_s.
	\]
\end{lemma}

\begin{lemma}\label{space}
For the auxiliary space \(X\), we have
	\begin{equation}\label{spaceX}
		X = Y \cap W^{1,p_1}(\mathbb{R}^N),
	\end{equation}
and the norms \(\|\cdot\|_X\) and
	\(\|u\|_Y+\|u\|_{W^{1,p_1}}\)
induce the same notion of convergence on this space.
\end{lemma}

\begin{proof}
Write
	\[
	Z=Y\cap W^{1,p_1}(\mathbb R^N),
	\]
Equip \(Z\) with the norm
	\[
	\|u\|_Z=\|u\|_Y+\|u\|_{W^{1,p_1}}.
	\]
Then \(Z\) is a Banach space.
	
We first prove that, on \(C_c^\infty(\mathbb R^N)\), the norm \(\|\cdot\|_X\)
is equivalent to \(\|\cdot\|_Z\). By definition,
	\[
	\|u\|_Y\le \|u\|_X,
	\qquad
	\|Du\|_{p_1}\le \|u\|_X.
	\]
It remains to control the function term \(\|u\|_{p_1}\) in the norm \(W^{1,p_1}\). Applying Lemma
\ref{lem:Gn} with \(d=N,\ s=p_1,\ q+1=2\) gives a constant \(C>0\) such that
	\[
	\|u\|_{\infty}
	\le
	C\|Du\|_{p_1}^{\theta}\|u\|_2^{1-\theta},
	\qquad
	\theta=
	\frac{Np_1}{Np_1+2p_1-2N}.
	\]
Therefore,
	\[
	\|u\|_{p_1}^{p_1}
	=
	\int_{\mathbb R^N}|u|^{p_1}\,dx
	\le
	\|u\|_\infty^{p_1-2}\|u\|_2^2.
	\]
Combining this with \(Y\hookrightarrow H^1(\mathbb R^N)\) yields
	\[
	\|u\|_{p_1}^{p_1}
	\le
	C
	\|Du\|_{p_1}^{\theta(p_1-2)}
	\|u\|_Y^{p_1-\theta(p_1-2)}.
	\]
Young's inequality then gives
	\[
	\|u\|_{p_1}
	\le
	C'\bigl(\|Du\|_{p_1}+\|u\|_Y\bigr).
	\]
Consequently,
	\(\|u\|_{W^{1,p_1}}
	\le
	C''\|u\|_X.\)
and hence
	\(\|u\|_Z\le C\|u\|_X.\)
Conversely, the definition immediately gives
	\(\|u\|_X\le \|u\|_Z.\)
Thus \(\|\cdot\|_X\) and \(\|\cdot\|_Z\) are equivalent on \(C_c^\infty(\mathbb R^N)\).
	
We next show that \(C_c^\infty(\mathbb R^N)\) is dense in \(Z\). Since
\(V\in C(\mathbb R^N)\) and \(V(x)\ge V_0>0\), for every
\(u\in Z\) one may first choose cutoffs \(u_R=\eta_R u\) such that \(u_R\to u\) both in \(Y\) and in \(W^{1,p_1}(\mathbb R^N)\).
Indeed,
	\[
	\int_{\mathbb R^N}V(x)u^2\,dx<+\infty,
	\qquad
	u\in W^{1,p_1}(\mathbb R^N),
	\]
One then mollifies \(u_R\) on a bounded domain; see Evans \cite{RN88}. Since
\(V\) is bounded on bounded sets, functions in \(C_c^\infty(\mathbb R^N)\)
approximate \(u_R\) in the \(\|\cdot\|_Y+\|\cdot\|_{W^{1,p_1}}\) norm. Therefore,
\(C_c^\infty(\mathbb R^N)\) is dense in \(Z\).
	
The equivalence of \(\|\cdot\|_X\) and \(\|\cdot\|_Z\) on
\(C_c^\infty(\mathbb R^N)\) now implies
	\[
	X
	=
	\overline{C_c^\infty(\mathbb R^N)}^{\|\cdot\|_X}
	=
	\overline{C_c^\infty(\mathbb R^N)}^{\|\cdot\|_Z}
	=
	Z.
	\]
that is,
	\(X=Y\cap W^{1,p_1}(\mathbb R^N).\)
\end{proof}

To identify elements of \(X\) with locally continuous functions, we recall Morrey's embedding on bounded domains \cite{RN88} and its consequence for the present space \(X\).

\begin{lemma}\label{lem:morrey-local}
Let \(q>N\) and set \(\alpha=1-\frac{N}{q}\in(0,1).\)
For every \(R>0\), the continuous embedding
	\begin{equation}\label{morrey}
		W^{1,q}(B_R)\hookrightarrow C^{0,\alpha}(\overline{B_R}).
	\end{equation}
holds. In particular, each \(u\in W^{1,q}(B_R)\) has a representative of \(u\)
\(u^*\) such that \(u=u^*\in C^{0,\alpha}(\overline{B_R})\) almost everywhere in \(B_R\).
Moreover, there is a constant \(C_R>0\) such that
	\[
	\|u^*\|_{C^{0,\alpha}(\overline{B_R})}
	\le
	C_R\|u\|_{W^{1,q}(B_R)}.
	\]
In particular, for every \(x,y\in\overline{B_R}\),
	\[
	|u^*(x)-u^*(y)|
	\le
	C_R\|u\|_{W^{1,q}(B_R)}|x-y|^\alpha.
	\]
\end{lemma}

We apply this to \(X\). Our choice of \(p_1\) gives \(p_1>N,\)
and Lemma \ref{space} gives \(X\hookrightarrow W^{1,p_1}_{loc}(\mathbb R^N).\)
Thus, for every \(R>0\) and every \(u\in X\),
\(u\in W^{1,p_1}(B_R).\)
Lemma \ref{lem:morrey-local} shows that \(u\) has a representative on \(B_R\),
still denoted by \(u\), satisfying \(u\in C^{0,\alpha}(\overline{B_R}),\)
Therefore,
\begin{equation}\label{morreyembed}
	X\hookrightarrow C^{0,\alpha}_{loc}(\mathbb R^N),
	\qquad
	X\hookrightarrow L^\infty_{loc}(\mathbb R^N).
\end{equation}
Whenever pointwise values of \(u\in X\) are used below, we take this
\(C^{0,\alpha}_{loc}(\mathbb R^N)\) continuous representative.

In the space \(X\), the effective domain of \(\Psi\) is naturally given by
\[
D_X(\Psi)
=
\{u\in X:\Psi(u)<+\infty\}
=
\{u\in X:\|Du\|_\infty\le 1\}.
\]
We shall abbreviate \(D_X(\Psi)\) as \(D(\Psi)\).
Under the present assumptions,
\begin{equation}
	D(\Psi)=D_Y(\Psi).
\end{equation}
Indeed, \(D(\Psi)\subset D_Y(\Psi)\) is immediate. Conversely, if
\(u\in D_Y(\Psi)\), then \(u\in Y\) and \(\|Du\|_\infty\le1\). Since
\(p_1\ge2\), we have
\[
|Du|^{p_1}
=
|Du|^{p_1-2}|Du|^2
\le
|Du|^2
\]
almost everywhere. Hence \(Du\in L^{p_1}(\mathbb R^N).\)
Since \(u\in Y\hookrightarrow H^1(\mathbb R^N)\) and
\(Du\in L^{p_1}(\mathbb R^N)\), Lemma \ref{lem:Gn} implies
\(u\in L^\infty(\mathbb R^N)\). Therefore,
\[
\|u\|_{p_1}^{p_1}
\le
\|u\|_\infty^{p_1-2}\|u\|_2^2<+\infty.
\]
Therefore,
\(u\in W^{1,p_1}(\mathbb R^N).\)
Combining this with \(u\in Y\) yields
\[
u\in Y\cap W^{1,p_1}(\mathbb R^N)=X.
\]
and hence \(u\in D_X(\Psi)\). This proves the claim.

On \(D(\Psi)\), convergence in the two spaces is equivalent, and critical points obtained in the auxiliary space are also critical in the natural space.
\begin{lemma}\label{lem:natural-auxiliary-domain}
Suppose that there exists a constant \(C>0\) such that
	\(\|Du_j\|_\infty\le C, \|Du\|_\infty\le C,\)
Then \(u_j\to u\) in \(X\) if and only if \(u_j\to u\) in \(Y\).
Moreover, if \(u\in D(\Psi)=D_Y(\Psi)\) is a critical point of \(I_\lambda\) in \(X\), then \(u\) is also a critical point of \(I_\lambda\) in the natural space \(Y\).
\end{lemma}

\begin{proof}
First, since \(X\subset Y\) and \(\|u\|_Y\le\|u\|_X\), if
\(u_j\to u\) in \(X\), then \(u_j\to u\) in \(Y\).
	
Conversely, suppose that \(u_j\to u\) in \(Y\). The hypothesis
	\[
	\|Du_j\|_\infty\le C,
	\qquad
	\|Du\|_\infty\le C,
	\]
gives
	\[
	\|D(u_j-u)\|_\infty\le2C.
	\]
Since \(p_1\ge2\),
	\[
	|D(u_j-u)|^{p_1}
	=
	|D(u_j-u)|^{p_1-2}|D(u_j-u)|^2
	\le
	(2C)^{p_1-2}|D(u_j-u)|^2.
	\]
Consequently,
	\[
	\|D(u_j-u)\|_{p_1}^{p_1}
	\le
	(2C)^{p_1-2}\|D(u_j-u)\|_2^2\to0.
	\]
and \(u_j\to u\) in \(Y\), we obtain \(\|u_j-u\|_Y\to0.\)
By the definition of \(\|\cdot\|_X\),
	\(\|u_j-u\|_X\to0.\)
Thus \(u_j\to u\) in \(X\).
	
If \(u\in D(\Psi)\) is a critical point of \(I_\lambda\) in \(X\), then, for every \(v\in D(\Psi)\),
	\[
	\lambda\bigl(\Psi(v)-\Psi(u)\bigr)
	-
	\Phi'(u)(v-u)
	\ge0.
	\]
Because \(D(\Psi)=D_Y(\Psi)\), this inequality also holds for every
\(v\in D_Y(\Psi)\). Thus \(u\) is critical for \(I_\lambda\) in the natural space \(Y\).
\end{proof}

\begin{remark}
Consequently, convergence, compactness, and critical point statements obtained in \(X\) transfer to the corresponding natural-space statements in \(Y\).
\end{remark}
The nonsmooth critical points constructed in \(X\) must ultimately be interpreted as weak solutions of \eqref{biff}. We therefore specify the weak formulation.
\begin{definition}\label{def:ws}
A function \(u\in X\) is a weak solution of \eqref{biff} if
\(u\in \operatorname{Lip}(\mathbb R^N)\) and, almost everywhere in \(\mathbb R^N\),
\(|Du|\leq 1 \),
and
	\[
	\frac{1}{\sqrt{1-|Du|^2}}\in L^1_{\mathrm{loc}}(\mathbb R^N),
	\]
and, writing
	\begin{equation}\label{eq:rhou}
		\rho_u(x)=|u|^{p-2}u-V(x)u,
	\end{equation}
for every \(\eta\in \operatorname{Lip}_c(\mathbb R^N)
	\coloneq
	\operatorname{Lip}(\mathbb R^N)\cap C_c(\mathbb R^N),\)
the following integral identity holds:
	\begin{equation}\label{eq:ws}
		\int_{\mathbb{R}^N}
		\frac{Du \cdot D\eta}{\sqrt{1 - |Du|^2}}\,dx
		=
		\langle \rho_u, \eta \rangle,
	\end{equation}
Here \(\langle \rho_u, \eta \rangle\) is the duality pairing, satisfying
	\[
	\langle \rho_u, \eta \rangle
	=
	\int_{\mathbb R^N}
	\left(|u|^{p-2}u-V(x)u\right)\eta\,dx.
	\]
\end{definition}

We next discuss the two parts \(\Psi\) and \(\Phi\) of the energy. For \(\Psi\), we have
\begin{lemma}\label{lem:Dpsi-weakly-closed}
	\(D(\Psi)\)
is a convex, weakly closed subset of \(X\).
\end{lemma}

\begin{proof}
The convexity of \(D(\Psi)\) is immediate.
We prove weak closedness.
	
Suppose that, in \(X\),
	\(u_j\rightharpoonup u\), \(u_j\in D(\Psi)\). 
Since \(X\) embeds continuously into both \(Y\) and \(W^{1,p_1}(\mathbb R^N)\),
we have \(u_j\rightharpoonup u\) in \(Y\) and \(	Du_j\rightharpoonup Du\) in \(W^{1,p_1}(\mathbb R^N)\).
The almost-everywhere bound \(|Du_j|\le1\) implies
	\[
	Du_j\in
	\mathcal K
	\coloneq
	\left\{
	F\in L^{p_1}(\mathbb R^N;\mathbb R^N):
	\|F\|_{\infty} \le 1
	\right\}.
	\]
The set \(\mathcal K\) is closed and convex in \(L^{p_1}(\mathbb R^N;\mathbb R^N)\), hence weakly closed. Therefore,
	\(Du\in\mathcal K,\)
that is, \(|Du|\le1\) almost everywhere on \(\mathbb R^N\). Also \(u\in X\), so \(u\in D(\Psi).\)
Thus \(D(\Psi)\) is weakly closed in \(X\).
\end{proof}

\subsection{Properties of Palais--Smale sequences}
\begin{lemma}[Right continuity of the mountain pass level at \(\lambda=1\)]\label{lem:c0-right-continuity}
For the family \(I_\lambda(u)\),
let \(c_0(\lambda)\) be the mountain pass level associated with the fixed path class \(\Gamma_0\).
Then \(c_0\) is right-continuous at \(\lambda=1\), namely
	\[
	\lim_{\lambda\downarrow1}c_0(\lambda)=c_0(1).
	\]
\end{lemma}

\begin{proof}
Since \(\Psi\ge0\), for \(\lambda\ge1\) and every \(u\in X\),
	\[
	I_\lambda(u)=I_1(u)+(\lambda-1)\Psi(u)\ge I_1(u).
	\]
As \(\Gamma_0\) does not depend on \(\lambda\), it follows that \(c_0(\lambda)\ge c_0(1)\).
It remains to prove
	\[
	\limsup_{\lambda\downarrow1}c_0(\lambda)\le c_0(1).
	\]
	
Fix \(\varepsilon>0\). By the definition of \(c_0(1)\), there is a path
\(\gamma_\varepsilon\in\Gamma_0\) such that
	\[
	\sup_{t\in[0,1]}I_1(\gamma_\varepsilon(t))
	\le c_0(1)+\varepsilon.
	\]
Since \(\gamma_\varepsilon([0,1])\) is compact in \(X\) and \(\Phi\) is continuous on \(X\),
	\(\sup_{t\in[0,1]}\Phi(\gamma_\varepsilon(t))<+\infty.\)
Using \(I_1(\gamma_\varepsilon(t))
	=
	\Psi(\gamma_\varepsilon(t))-\Phi(\gamma_\varepsilon(t))\)
gives
	\[
	M_\varepsilon
	\coloneq
	\sup_{t\in[0,1]}\Psi(\gamma_\varepsilon(t))<+\infty.
	\]
we obtain, as \(\lambda\ge1\),
	\[
	\begin{aligned}
		c_0(\lambda)
		&\le
		\sup_{t\in[0,1]}I_\lambda(\gamma_\varepsilon(t))\\
		&=
		\sup_{t\in[0,1]}
		\left[
		I_1(\gamma_\varepsilon(t))
		+
		(\lambda-1)\Psi(\gamma_\varepsilon(t))
		\right]\\
		&\le
		c_0(1)+\varepsilon+(\lambda-1)M_\varepsilon.
	\end{aligned}
	\]
Letting \(\lambda\downarrow1\) and then \(\varepsilon\to0\),
gives
	\[
	\limsup_{\lambda\downarrow1}c_0(\lambda)\le c_0(1).
	\]
and combining this with \(c_0(\lambda)\ge c_0(1)\), proves the assertion.
\end{proof}

\subsection{Auxiliary scaling and geometric formulas}\label{compute}
We first record the basic inequality used in scaling comparisons \cite{WOS:001150613900001}:
\begin{lemma}\label{lem:BI-basic}
Let \(0\le t\le 1\). Then
	\begin{equation}\label{foundation}
		\frac{t^2}{2}
		\le
		1-\sqrt{1-t^2}
		\le
		t^2.
	\end{equation}
Moreover, the function \(	t\mapsto 1-\sqrt{1-t^2}\)
admits the expansion, for \([0,1)\),
	\begin{equation}\label{zk}
		1-\sqrt{1-t^2}
		=
		\sum_{k=1}^{\infty} b_k t^{2k},
		\qquad b_k>0.
	\end{equation}
For points satisfying \(|Du|<1\) and \(|D\psi|\le1\), one also has
	\begin{equation}
		\frac{Du\cdot(Du-D\psi)}{\sqrt{1-|Du|^2}}
		\ge
		\sqrt{1-|D\psi|^2}-\sqrt{1-|Du|^2}.
	\end{equation}
\end{lemma}

For the regularity argument, we recall the elementary geometric formulas for spacelike graphs in Lorentz--Minkowski space.
The induced metric, second fundamental form, and mean curvature depend only on
\(u(x)\), not on the particular source in the equation.
We therefore use the calculations in \cite{WOS:001150613900001} and specify the relation between the mean curvature \(H\) and the source in our non-autonomous problem.

\begin{definition}\label{picture}
Let \(u\in C^2(\mathbb R^N)\) satisfy \(|Du|<1\) on \(\mathbb R^N\),
and consider Lorentz--Minkowski space
	\(\mathbb L^{N+1}=
\mathbb R\times\mathbb R^N\). Define the graph
	\(M=\{(u(x),x):x\in\mathbb R^N\}.\)
and its parametrization
	\begin{equation}\label{F}
		F(x)=(u(x),x).
	\end{equation}
\end{definition}

\begin{proposition}
Under the assumptions of Definition \ref{picture},
the mean curvature satisfies
	\begin{equation}\label{eq:mean-curvature-geometric}
		H
		=
		\operatorname{div}
		\left(
		\frac{Du}{\sqrt{1-|Du|^2}}
		\right).
	\end{equation}
If, moreover, \(u\) is a classical solution of \eqref{hceq}, then
the identity \eqref{H} holds, namely
	\begin{equation}\label{H2}
		H_\lambda
		=
		\operatorname{div}\left(\frac{Du}{\sqrt{1-|Du|^2}}\right)
		=
		V(x)u-\lambda^{-1}|u|^{p-2}u.
	\end{equation}
\end{proposition}
\begin{proof}
Choose coordinates \(x=(x^1,\ldots,x^N)\) and write
	\[
	u_i=\frac{\partial u}{\partial x^i},
	\qquad
	u_{ij}=\frac{\partial^2u}{\partial x^i\partial x^j}.
	\]
The tangent vectors are \(F_i
	=
	F_*\partial_i
	=
	(u_i,e_i).\)
With respect to the Minkowski metric \(\langle (t,x),(s,y)\rangle_{\mathbb L}
	=
	-ts+x\cdot y\)
the induced metric is
	\begin{equation}\label{gij}
		g_{ij}
		=
		\langle F_i,F_j\rangle_{\mathbb L}
		=
		\delta_{ij}-u_i u_j.
	\end{equation}
Set
	\begin{equation}\label{w}
		w
		=
		\frac1{\sqrt{1-|Du|^2}},
	\end{equation}
Then the inverse metric and volume form are
	\[
	g^{ij}
	=
	\delta_{ij}+w^2u_i u_j,
	\qquad
	dx_g=w^{-1}dx.
	\]
The second fundamental form \(A=(h_{ij})\) satisfies
	\(h_{ij}
	=
	wu_{ij}.\)
and its norm is
	\[
	|A|_g^2
	=
	w^2|D^2u|^2
	+
	2w^4|D^2u\,Du|^2
	+
	w^6\bigl(D^2u[Du,Du]\bigr)^2.
	\]
	
With these conventions, the Lorentzian mean curvature of the graph \(u(x)\) is the trace of its second fundamental form with respect to the induced metric, namely \(H
	=
	g^{ij}h_{ij}.\)
This gives
	\begin{equation}
		\begin{aligned}
			\operatorname{div}\left(\frac{Du}{\sqrt{1-|Du|^2}}\right)
			&=\sum_{i=1}^N \partial_i(wu_i)  \\
			&= w\Delta u+w^3D^2u[Du,Du]\\
			&=H.
		\end{aligned}
	\end{equation}
Substituting equation \eqref{hceq} into the left-hand side yields the formula \eqref{H}.
\end{proof}

\section{Construction of parameter-dependent critical points}
We now construct critical points of the parameter-dependent functional and verify the required hypotheses.
In the nonsmooth framework of Section 1, we first verify on the auxiliary space \(X\) that
the functional \(I_\lambda=\lambda\Psi-\Phi\)
satisfies \((\Psi_1)\), \((\Psi_2)\), \((\Psi_3)\), \((\Phi_1)\), and \((\Phi_2)\).
To apply the monotonicity trick and the multiplicity theorem, we also prove, for \(\lambda\) near \(1\),
the family \(I_\lambda\) has uniform mountain pass and symmetric mountain pass geometry, together with monotonicity and finiteness of the minimax levels.

The ordinary mountain pass construction gives bounded PS sequences at \(c_0(\lambda)\), which converge strongly for each fixed parameter to a critical point.
For \(c_k(\lambda)\), we verify only the geometric and monotonicity hypotheses needed for Lemma \ref{th1.7},
without constructing a separate parameter-dependent critical-point sequence for each \(k\ge1\).

\subsection{Verification of the functional assumptions}
We first verify that \(X\), \(\Phi\), and \(\Psi\)
satisfy the assumptions used below.
\begin{lemma}\label{space-properties}
Let \(p_1>\max\{N, p\}\) and let \(X\) be defined in \eqref{fzr}. Then \(X\) is a reflexive Banach space.
For every \(2\le q\le 2^*\), the embedding
	\(X\hookrightarrow L^q(\mathbb{R}^N).\)
is continuous, and, if \(2\le q<2^*\), the embedding
\(X\hookrightarrow L^q(\mathbb R^N)\) is compact.
Finally,
	\(X\hookrightarrow C_0(\mathbb{R}^N),\)
where
	\[
	C_0(\mathbb{R}^N)
	=
	\left\{
	u\in C(\mathbb{R}^N):
	\lim_{|x|\to\infty}u(x)=0
	\right\}.
	\]
\end{lemma}

\begin{proof}
The Hilbert space \(Y=H_V^1(\mathbb R^N)\) is reflexive.
Since \(1<p_1<+\infty\), the space \(W^{1,p_1}(\mathbb R^N)\) is reflexive as well.
By Lemma \ref{space}, \(X=Y\cap W^{1,p_1}(\mathbb R^N),\)
and \(\|\cdot\|_X\) is equivalent to the natural intersection norm.
Hence \(X\) is a reflexive Banach space.
	
The assumption \(V(x)\ge V_0>0\) implies
for some constant \(C>0\).
	\[
	\|u\|_{H^1(\mathbb R^N)}
	\le
	C\|u\|_Y
	\le
	C\|u\|_X.
	\]
Combining this with the Sobolev embedding
	\[
	H^1(\mathbb{R}^N)\hookrightarrow L^q(\mathbb{R}^N),
	\qquad 2\le q\le 2^*,
	\]
gives, for every \(u\in X\), a constant \(C_q>0\) such that
	\(\|u\|_q\le C_q\|u\|_X.\)
Thus \(X\hookrightarrow L^q(\mathbb R^N)\) continuously.
	
For compactness when \(2\le q<2^*\), let \(\{u_j\}\) be bounded in \(X\).
Then \(\{u_j\}\) is bounded in \(Y\) and, in particular, in \(H^1(\mathbb R^N)\).
Local Rellich compactness gives, for each \(R>0\), a subsequence such that
\(u_j\to u\) in \(L^q(B_R(0))\).
	
On the other hand, \(V(x)\to+\infty\) implies that, for any \(\varepsilon>0\), one can choose \(R>0\) such that
	\(\inf_{|x|>R}V(x)\)
is sufficiently large. Boundedness of \(\{u_j\}\) in \(Y\) then yields
	\begin{equation}\label{yuxiang}
		\int_{|x|>R}|u_j|^2\,dx
		\le
		\frac1{\inf_{|x|>R}V(x)}
		\int_{|x|>R}V(x)u_j^2\,dx
		\le
		\varepsilon
	\end{equation}
uniformly for all \(j\). Together with \(H^1\hookrightarrow L^{2^*}\), there is a constant \(C_*>0\) such that, for
\(\forall j\), one has \(\|u_j\|_{2^*}\le C_*,\)
For \(2\le q<2^*\), choose \(\theta\in(0,1]\) satisfying
	\[
	\frac1q=\frac{\theta}{2}+\frac{1-\theta}{2^*},
	\]
By H{\"o}lder interpolation, for every measurable set \(E\subset\mathbb R^N\),
	\[
	\|u_j\|_{L^q(E)}\le \|u_j\|_{L^2(E)}^{\theta}\|u_j\|_{L^{2^*}(E)}^{1-\theta}.
	\]
Taking \(E=\{x:|x|\ge R\}\) and using \eqref{yuxiang} gives
	\[
	\|u_j\|_{L^q(|x|\ge R)}
	\le \bigl(\|u_j\|_{L^2(|x|\ge R)}\bigr)^{\theta}\bigl(\|u_j\|_{2^*}\bigr)^{1-\theta}
	\le \varepsilon^{\theta}C_*^{1-\theta}.
	\]
Consequently,
	\(\int_{|x|>R}|u_j|^q\,dx\)
uniform smallness for all \(j\). Combined with local strong \(L^q\) convergence, this gives, in \(L^q(\mathbb R^N)\),
	\(u_j\to u\), 
Therefore \(X\hookrightarrow L^q(\mathbb R^N)\) compactly.
	
Moreover, \(X\hookrightarrow W^{1,p_1}(\mathbb R^N)\) and Lemma \ref{lem:morrey-local} give
	\[
	X\hookrightarrow C_{\mathrm{loc}}^{0,\alpha}(\mathbb{R}^N),
	\qquad
	\alpha=1-\frac{N}{p_1}.
	\]
Lemma \ref{lem:Gn} also yields a constant \(C>0\) such that
	\[
	\|u\|_{\infty}
	\le
	C\|Du\|_{p_1}^{\theta}\|u\|_2^{1-\theta}
	\le
	C\|u\|_X.
	\]
Since \(X\) is the closure of \(C_c^\infty(\mathbb R^N)\) in the \(\|\cdot\|_X\) norm,
while \(C_c^\infty(\mathbb R^N)\subset C_0(\mathbb R^N)\) and convergence in \(X\) implies uniform convergence,
we conclude that \(X\hookrightarrow C_0(\mathbb R^N)\).
\end{proof}

\begin{lemma}\label{prop:Phi-Psi-conditions}
The functional \(\Phi: X\to\mathbb{R}\) satisfies \((\Phi_1)\),
and \(\Psi: X\to[0,+\infty]\) satisfies \((\Psi_1)\), \((\Psi_2)\), and \((\Psi_3)\).
\end{lemma}

\begin{proof}
We first verify \((\Phi_1)\).
	
By Lemma \ref{space-properties},
\(X\hookrightarrow L^p(\mathbb R^N)\) continuously.
Thus, for every \(u\in X\),
	\[
	|\Phi(u)|
	\le
	\frac1p\|u\|_p^p
	\le
	C\|u\|_X^p<+\infty.
	\]
so \(\Phi\) is well defined on \(X\).
	
H{\"o}lder's inequality gives, for every \(u,v\in X\),
	\[
	\begin{aligned}
		|\Phi'(u)v|
		&=
		\left|
		\int_{\mathbb R^N}|u|^{p-2}uv\,dx
		\right| \\
		&\le
		\int_{\mathbb R^N}|u|^{p-1}|v|\,dx \\
		&\le
		\|u\|_p^{p-1}\|v\|_p \\
		&\le
		C\|u\|_X^{p-1}\|v\|_X.
	\end{aligned}
	\]
Therefore \(\Phi'(u)\in X^*\).
	
Suppose that \(u_n\to u\) in \(X\).
Then \(u_n\to u\) in \(L^p(\mathbb R^N)\), and hence, in
\(L^{\frac{p}{p-1}}(\mathbb R^N)\),
	\[
	|u_n|^{p-2}u_n\to |u|^{p-2}u.
	\]
Therefore,
	\[
	\begin{aligned}
		\|\Phi'(u_n)-\Phi'(u)\|_{X^*}
		&=
		\sup_{\|v\|_X\le1}
		\left|
		\int_{\mathbb R^N}
		\bigl(|u_n|^{p-2}u_n-|u|^{p-2}u\bigr)v\,dx
		\right|  \\
		&\le
		\bigl\|
		|u_n|^{p-2}u_n-|u|^{p-2}u
		\bigr\|_{\frac{p}{p-1}}
		\sup_{\|v\|_X\le1}\|v\|_p
		\to0.
	\end{aligned}
	\]
Thus \(\Phi\in C^1(X,\mathbb R)\), proving \((\Phi_1)\).
	
We turn to the properties of \(\Psi\).
Clearly \(\Psi(0)=0\). Its convexity follows from Lemma 2.17 of \cite{Byeon2026} and the convexity of
	\(u\mapsto \frac12\int_{\mathbb R^N}V(x)u^2\,dx\)
the quadratic potential term.
	
For lower semicontinuity of \(\Psi\), suppose that \(u_j\to u\) in \(X\).
If
	\(\liminf_{j\to\infty}\Psi(u_j)=+\infty,\)
the conclusion is immediate. Otherwise,
	\[
	\liminf_{j\to\infty}\Psi(u_j)<+\infty,
	\]
choose a subsequence, still denoted by \(\{u_j\}\), such that
	\[
	\Psi(u_j)\to \liminf_{j\to\infty}\Psi(u_j)<+\infty.
	\]
Then \(|Du_j(x)|\le1\) almost everywhere in \(\mathbb R^N\).
Since \(u_j\to u\) in \(X\), Lemma \ref{space} implies
\(u_j\to u\) in \(Y\) and, in
\(L^{p_1}(\mathbb R^N;\mathbb R^N)\),
	\(Du_j\to Du\).
a further subsequence, \(Du_j(x)\to Du(x)\) almost everywhere on \(\mathbb R^N\).
Consequently \(|Du|\le1\) almost everywhere, that is, \(u\in D(\Psi)\).
	
Also \(u_j\to u\) in \(Y\) implies
	\[
	\frac12\int_{\mathbb R^N}V(x)u_j^2\,dx
	\to
	\frac12\int_{\mathbb R^N}V(x)u^2\,dx.
	\]
Fatou's lemma gives
	\[
	\int_{\mathbb R^N}
	\left(1-\sqrt{1-|Du|^2}\right)\,dx
	\le
	\liminf_{j\to\infty}
	\int_{\mathbb R^N}
	\left(1-\sqrt{1-|Du_j|^2}\right)\,dx.
	\]
Therefore,
	\[
	\Psi(u)\le \liminf_{j\to\infty}\Psi(u_j).
	\]
Hence \(\Psi\) is lower semicontinuous on \(X\).
This proves \((\Psi_1)\).
	
For \((\Psi_2)\), if \(u\in D(\Psi)\), then \(\|Du\|_\infty\le1\).
Using \eqref{foundation}, we obtain
	\[
	\int_{\mathbb R^N}
	\left(1-\sqrt{1-|Du|^2}\right)\,dx
	\ge
	\frac12\|Du\|_2^2.
	\]
Consequently,
	\begin{equation}\label{psiu}
		\Psi(u)
		\ge
		\frac12
		\left(
		\|Du\|_2^2
		+
		\int_{\mathbb R^N}V(x)u^2\,dx
		\right)
		=
		\frac12\|u\|_Y^2.
	\end{equation}
Since \(\|Du\|_\infty\le1\) when \(u\in D(\Psi)\), and \(p_1\ge2\),
	\[
	\|Du\|_{p_1}^{p_1}
	\le
	\|Du\|_2^2
	\le
	\|u\|_Y^2.
	\]
The condition \(\|u\|_X^2=\|u\|_Y^2+\|Du\|_{p_1}^2\) implies
If \(\|u\|_X\to\infty\), then necessarily \(\|u\|_Y\to\infty\).
Together with \eqref{psiu}, this gives \(\Psi(u)\to+\infty.\)
and proves \((\Psi_2)\).
	
Finally, we verify \((\Psi_3)\).
Suppose that \(u_j\rightharpoonup u\) in \(X\) and \(\Psi(u_j)\to\Psi(u)\).
We may assume \(u_j,u\in D(\Psi)\). The expansion in Lemma \ref{lem:BI-basic} and monotone convergence give
	\[
	\int_{\mathbb R^N}
	\left(1-\sqrt{1-|Du|^2}\right)\,dx
	=
	\sum_{k=1}^{\infty}
	b_k\|Du\|_{2k}^{2k}.
	\]
Consequently,
	\begin{equation}\label{split}
		\Psi(u)
		=
		\frac12
		\int_{\mathbb R^N}V(x)u^2\,dx
		+
		\sum_{k=1}^{\infty}b_k\|Du\|_{2k}^{2k}.
	\end{equation}
First, for every fixed \(k\in\mathbb N\), we claim that \(Du_j\rightharpoonup Du\) in \(L^{2k}(\mathbb R^N;\mathbb R^N)\).
From \(u_j\rightharpoonup u\), we already know
\(Du_j\rightharpoonup Du\) in \(L^2(\mathbb R^N;\mathbb R^N)\).
On the other hand, \(|Du_j|\le 1\) implies that, for every \(r\ge 2\),
\[
\|Du_j\|_r^r
=
\int_{\mathbb R^N}|Du_j|^r\,dx
\le
\int_{\mathbb R^N}|Du_j|^2\,dx.
\]
Thus \(\{Du_j\}\) is bounded in \(L^r(\mathbb R^N;\mathbb R^N)\). For any \(\varphi\in L^{r'}(\mathbb R^N;\mathbb R^N)\), choose
\(\varphi_m\in C_c^\infty(\mathbb R^N;\mathbb R^N)\) such that \(\varphi_m\) converges to \(\varphi\) in \(L^{r'}(\mathbb R^N;\mathbb R^N)\).
Since \(\varphi_m\in L^2(\mathbb R^N;\mathbb R^N)\),
\[
\int_{\mathbb R^N}(Du_j-Du)\cdot \varphi_m\,dx\to 0.
\]
and since \(\{Du_j-Du\}\) is uniformly bounded in \(L^r(\mathbb R^N;\mathbb R^N)\),
Therefore,
\[
\left|
\int_{\mathbb R^N}(Du_j-Du)\cdot(\varphi-\varphi_m)\,dx
\right|
\le
C\|\varphi-\varphi_m\|_{r'}.
\]
letting first \(j\to\infty\) and then \(m\to\infty\) gives
\[
\int_{\mathbb R^N}(Du_j-Du)\cdot \varphi\,dx\to 0.
\]
Therefore, for every \(r\ge 2\), \(Du_j\rightharpoonup Du\) in \(L^r(\mathbb R^N;\mathbb R^N)\).

Since \(u_j,u\in D(\Psi)\) and the \(L^q(1<q<\infty)\) norm is weakly lower semicontinuous, while
in \(Y\),
\(u_j\rightharpoonup u\), we have
	\[
	\int_{\mathbb R^N}V(x)u^2\,dx
	\le
	\liminf_{j\to\infty}
	\int_{\mathbb R^N}V(x)u_j^2\,dx,
	\]
and
	\[
	\|Du\|_{2k}^{2k}
	\le
	\liminf_{j\to\infty}\|Du_j\|_{2k}^{2k},
	\qquad
	\forall k\ge1.
	\]
	
For every \(m\in\mathbb N\), weak lower semicontinuity also gives
	\[
	\frac12\int_{\mathbb R^N}V(x)u^2\,dx
	+
	\sum_{k=1}^{m}b_k\|Du\|_{2k}^{2k}
	\le
	\liminf_{j\to\infty}
	\left[
	\frac12\int_{\mathbb R^N}V(x)u_j^2\,dx
	+
	\sum_{k=1}^{m}b_k\|Du_j\|_{2k}^{2k}
	\right].
	\]
If any of these inequalities were strict, the condition \(b_k>0\) would imply
	\[
	\Psi(u)<\liminf_{j\to\infty}\Psi(u_j)=\Psi(u),
	\]
a contradiction. Hence
	\[
	\int_{\mathbb R^N}V(x)u_j^2\,dx
	\to
	\int_{\mathbb R^N}V(x)u^2\,dx,
	\qquad
	\|Du_j\|_2\to\|Du\|_2.
	\]
It follows that, in \(L^2(\mathbb R^N;\mathbb R^N)\),
	\[
	Du_j\to Du,
	\qquad
	\sqrt{V}u_j\to \sqrt{V}u,
	\]
Thus \(u_j\to u\) in \(Y\).
Using \(u_j,u\in D(\Psi)\), we obtain
	\[
	\|D(u_j-u)\|_{p_1}^{p_1}
	\le
	2^{p_1-2}
	\|D(u_j-u)\|_2^2
	\to 0.
	\]
and consequently, in \(L^{p_1}(\mathbb R^N;\mathbb R^N)\),
	\(Du_j\to Du\). 
The definition of \(\|\cdot\|_X\) therefore gives
\(u_j\to u\) in \(X\).
This proves \((\Psi_3)\).
\end{proof}

\begin{lemma}\label{lem:Phi2-condition}
The functional \(\Phi\) satisfies \((\Phi_2)\).
\end{lemma}

\begin{proof}
Let \(\{u_j\}\) be a bounded sequence in \(X\). Since \(X\) is reflexive,
there exist \(u\in X\) and a subsequence, still denoted by \(\{u_j\}\), with \(u_j\rightharpoonup u\) in \(X\).
By Lemma \ref{space-properties}, after taking \(p'=\frac{p}{p-1}\),
we have, in \(L^{p'}(\mathbb R^N)\), \(|u_j|^{p-2}u_j
	\to
	|u|^{p-2}u\). 
	
For every \(v\in X\), the condition \(X\hookrightarrow L^p(\mathbb R^N)\)
and H{\"o}lder's inequality give
	\[
	\int_{\mathbb R^N}
	|u_j|^{p-2}u_jv\,dx
	\to
	\int_{\mathbb R^N}
	|u|^{p-2}uv\,dx.
	\]
On the other hand, \(u_j\to u\) in \(L^p(\mathbb R^N)\) implies
	\[
		\Phi'(u_j)(v-u_j)=
		\int_{\mathbb R^N}
		|u_j|^{p-2}u_j(v-u_j)\,dx
		\to
		\int_{\mathbb R^N}
		|u|^{p-2}u(v-u)\,dx
		=
		\Phi'(u)(v-u).
	\]
In particular,
	\[
	\liminf_{j\to\infty}
	\Phi'(u_j)(v-u_j)
	\ge
	\Phi'(u)(v-u).
	\]
Thus \(\Phi\) satisfies \((\Phi_2)\).
\end{proof}

We also verify evenness, condition (E).
\begin{lemma}\label{lem:even-condition}
Both \(\Psi\) and \(\Phi\) are even; hence the family
	\(I_\lambda\)
satisfies (E), namely
	\[
	I_\lambda(-u)=I_\lambda(u),
	\qquad
	\forall u\in X,\quad \forall \lambda>0.
	\]
\end{lemma}

\begin{proof}
By the definition of \(D(\Psi)\), if \(u\in D(\Psi)\),
then \(-u\in D(\Psi)\).
If \(u\notin D(\Psi)\), this equivalence gives \(-u\notin D(\Psi)\) and \(\Psi(u)=\Psi(-u)=+\infty.\)
	
For \(\Psi\),
	\[
	\begin{aligned}
		\Psi(-u)
		&=
		\int_{\mathbb R^N}
		\left(
		1-\sqrt{1-|D(-u)|^2}
		\right)\,dx
		+
		\frac12
		\int_{\mathbb R^N}V(x)(-u)^2\,dx\\
		&=
		\int_{\mathbb R^N}
		\left(
		1-\sqrt{1-|Du|^2}
		\right)\,dx
		+
		\frac12
		\int_{\mathbb R^N}V(x)u^2\,dx\\
		&=
		\Psi(u).
	\end{aligned}
	\]
Therefore, for every \(u\in X\),
	\(\Psi(-u)=\Psi(u).\)
The definition of \(\Phi\) also gives
	\[
	\Phi(-u)
	=
	\frac1p
	\int_{\mathbb R^N}|-u|^p\,dx
	=
	\frac1p
	\int_{\mathbb R^N}|u|^p\,dx
	=
	\Phi(u).
	\]
Thus, for every \(\lambda>0\),
	\[
	I_\lambda(-u)
	=
	\lambda\Psi(-u)-\Phi(-u)
	=
	\lambda\Psi(u)-\Phi(u)
	=
	I_\lambda(u).
	\]
Thus \(I_\lambda\) satisfies the evenness condition \((E)\).
\end{proof}

\subsection{Mountain pass geometry and critical-point construction}
We first show that the family
\(I_\lambda=\lambda\Psi-\Phi\)
has uniform mountain pass geometry on \(\lambda\in[1-\varepsilon_0,1+\varepsilon_0]\).
The lower bound on the sphere will also be used to establish uniform symmetric mountain pass geometry.
Using monotonicity of \(c_0(\lambda)\) in \(\lambda\), choose differentiability points \(\lambda_n\downarrow1\) and construct bounded PS sequences for \(I_{\lambda_n}\). Their strong convergence in \(X\) gives critical points of \(I_{\lambda_n}\) at level \(c_0(\lambda_n)\).

We begin with the required geometry.
\begin{lemma}\label{MP-Q}
The family \(I_\lambda\) in Definition \ref{def:critical} satisfies (MP):
there exist \(\varepsilon_0\in(0,\frac{1}{2})\), \(\rho_0>0\), \(\alpha_0>0\), and \(e\in D(\Psi)\) such that \(\|e\|_X>\rho_0\) and
	\[
	\inf_{\|u\|_X=\rho_0} I_{1-\varepsilon_0}(u)\ge \alpha_0>0,
	\qquad
	I_{1+\varepsilon_0}(e)<0.
	\]
\end{lemma}
	\begin{proof}
For $u\in D(\Psi)$, write
		\[
		a(u) = \|u\|_Y^2 = \int_{\mathbb R^N} \bigl(|Du|^2 + V(x)u^2\bigr)\,dx.
		\]
By \eqref{psiu} and Lemmas \ref{space-properties} and \ref{prop:Phi-Psi-conditions},
		\[
		\Phi(u) = \frac1p\int_{\mathbb R^N}|u|^p\,dx \le C\|u\|_p^p \le C\|u\|_Y^p = C a(u)^{p/2}.
		\]
Therefore, for any fixed $\varepsilon_0\in(0,1/2)$,
		\[
		I_{1-\varepsilon_0}(u)
		= (1-\varepsilon_0)\Psi(u)-\Phi(u)
		\ge \frac{1-\varepsilon_0}{2}a(u)-C a(u)^{p/2}.
		\]
Since $p>2$, the right-hand side is positive for sufficiently small $a(u)>0$.
		
We next ensure the required smallness of $a(u)$ on the sphere $\|u\|_X=\rho_0$.
For \(a(u)\le \rho_0^2<1\), choosing $\rho_0>0$ sufficiently small gives
		\[
		C a(u)^{p/2} \le \frac{1-\varepsilon_0}{4}a(u)
		\]
whenever $\|u\|_X=\rho_0$ and $u\in D(\Psi)$. Consequently,
		\[
		I_{1-\varepsilon_0}(u) \ge \frac{1-\varepsilon_0}{4}a(u).
		\]
On the other hand,
		\[
		\rho_0^2 = \|u\|_Y^2 + \|Du\|_{p_1}^2
		\le a(u) + a(u)^{2/p_1}
		\le 2 a(u)^{2/p_1}
		\]
implies \(a(u) \ge 2^{-p_1/2}\rho_0^{p_1}.\)
Consequently,
		\[
		I_{1-\varepsilon_0}(u)
		\ge \frac{1-\varepsilon_0}{4} 2^{-p_1/2}\rho_0^{p_1}.
		\]
Set
		\[
		\alpha_0 = \frac{1-\varepsilon_0}{4}2^{-p_1/2}\rho_0^{p_1} > 0,
		\]
and hence
		\[
		\inf_{\|u\|_X=\rho_0} I_{1-\varepsilon_0}(u) \ge \alpha_0 > 0.
		\]
	
We now construct the endpoint \(e\).
We have
	\[
	I_{1+\varepsilon_0}(e)
	=
	(1+\varepsilon_0)\int_{\mathbb R^N}G(De)\,dx
	+
	\frac{1+\varepsilon_0}{2}\int_{\mathbb R^N}V(x)e^2\,dx
	-
	\frac1p\int_{\mathbb R^N}|e|^p\,dx.
	\]
Choose \(\zeta\in C_c^\infty(\mathbb R^N)\) such that \(0 \le \zeta \le 1,\), \(\zeta \equiv 1\) on \(B_1(0)\), and \(\operatorname{supp} \zeta \subset B_2(0).\)
Choose a constant \(a>0\) such that
	\begin{equation}\label{eq:a-zeta}
		a\|D\zeta\|_\infty<1.
	\end{equation}
Set \(u_*(y)=a\zeta(y)\).
For \(t>0\), define
	\begin{equation}\label{eq:et-scaling}
		e_t(x)=t\,u_*\left(\frac{x}{t}\right)
		=
		at\,\zeta\left(\frac{x}{t}\right).
	\end{equation}
Then
	\[
	De_t(x)
	=
	Du_*\left(\frac{x}{t}\right)
	=
	aD\zeta\left(\frac{x}{t}\right).
	\]
By \eqref{eq:a-zeta},
	\(\|De_t\|_\infty
	=
	a\|D\zeta\|_\infty<1.\)
Thus \(e_t\in D(\Psi).\) for every \(t>0\).
	
We estimate \(I_{1+\varepsilon_0}(e_t)\). With the change of variables \(y=x/t\),
	\[
	\begin{aligned}
		\int_{\mathbb R^N}G(De_t)\,dx
		&=
		t^N
		\int_{\mathbb R^N}G(Du_*(y))\,dy\\
		&=
		t^N
		\int_{\mathbb R^N}G(aD\zeta(y))\,dy\\
		&= C_1t^N,
	\end{aligned}
	\]
Also,
	\[
	\begin{aligned}
		\int_{\mathbb R^N}|e_t|^p\,dx
		&=
		t^{N+p}\int_{\mathbb R^N}|u_*|^p\,dy\\
		&=
		a^pt^{N+p}\int_{\mathbb R^N}|\zeta|^p\,dy\\
		&= C_3t^{N+p},
	\end{aligned}
	\]
Furthermore,
	\[
	\int_{\mathbb R^N}V(x)e_t^2\,dx
	=
	a^2t^{N+2}
	\int_{\mathbb R^N}V(ty)\zeta^2(y)\,dy.
	\]
For fixed \(y\ne0\), set
	\begin{equation}\label{hys}
		h_y(s)=V(sy),\qquad s\ge1.
	\end{equation}
Using \eqref{V-derivative-growth}, we obtain
	\begin{equation}
		h_y'(s)
		=
		DV(sy)\cdot y
		=
		\frac{(sy)\cdot DV(sy)}{s}\le
		\frac{C_D}{s}h_y(s).
	\end{equation}
which can be written as
	\begin{equation}
		\frac{d}{ds}\log h_y(s)\le \frac{C_D}{s}.
	\end{equation}
Integrating over \(s\in[1,t]\) gives
	\begin{equation}\label{vtytcd}
		V(ty)=h_y(t)\le t^{C_D}h_y(1)=t^{C_D}V(y),
		\qquad t\ge1.
	\end{equation}
Since \(\zeta\) has compact support and \(V\in C^1(\mathbb R^N)\),
	\[
	\int_{\mathbb R^N}V(y)\zeta^2(y)\,dy<+\infty.
	\]
we obtain, for \(t\ge1\),
	\[
	\begin{aligned}
		\int_{\mathbb R^N}V(x)e_t^2\,dx
		&=
		a^2t^{N+2}
		\int_{\mathbb R^N}V(ty)\zeta^2(y)\,dy\\
		&\le
		a^2t^{N+2+C_D}
		\int_{\mathbb R^N}V(y)\zeta^2(y)\,dy\\
		&= C_2t^{N+2+C_D},
	\end{aligned}
	\]
	
Combining the preceding estimates yields
	\[
	\begin{aligned}
		I_{1+\varepsilon_0}(e_t)
		&=
		(1+\varepsilon_0)\Psi(e_t)-\Phi(e_t)\\
		&\le
		(1+\varepsilon_0)C_1t^N
		+
		\frac{1+\varepsilon_0}{2}C_2t^{N+2+C_D}
		-
		\frac{C_3}{p}t^{N+p}.
	\end{aligned}
	\]
Dividing the right-hand side by \(t^{N+p}\) gives
	\begin{equation}
		\frac{I_{1+\varepsilon_0}(e_t)}{t^{N+p}}
		\le
		(1+\varepsilon_0)C_1t^{-p}
		+
		\frac{1+\varepsilon_0}{2}C_2t^{C_D+2-p}
		-
		\frac{C_3}{p}.
	\end{equation}
Since \(C_D<p-2\), there exists \(T_1>0\) such that, for \(t\ge T_1\),
	\[
	\frac{I_{1+\varepsilon_0}(e_t)}{t^{N+p}}
	\le
	-\frac{C_3}{2p}<0.
	\]
Hence \(I_{1+\varepsilon_0}(e_t)<0\) whenever \(t\ge T_1\).
	
For the norm, observe that
	\[
	\|e_t\|_X^2
	\ge
	\|De_t\|_2^2
	=
	a^2t^N\|D\zeta\|_2^2.
	\]
There exists \(T_2>0\) such that, for \(t\ge T_2\),
	\(\|e_t\|_X>\rho_0.\)

Choose \(T>\max\{T_1,T_2,1\},\) and set \(e=e_T.\)
Then \(e\in D(\Psi),\)
and both \(\|e\|_X>\rho_0\) and \(I_{1+\varepsilon_0}(e)<0\) hold.
Thus \(e\) is the negative-energy endpoint required by the mountain pass geometry.
\end{proof}
\begin{remark}\label{ae}
By Definition \ref{def:mountain}, for every
\(\lambda\in[1-\varepsilon_0,1+\varepsilon_0]\), use the fixed endpoint \(e\) above to define
	\[
	\Gamma_0
	=
	\{\gamma\in C([0,1],X):\gamma(0)=0,\ \gamma(1)=e\}.
	\]
Since
	\[
	I_{1+\varepsilon_0}(e)<0
	\]
and \(\Psi(e)\ge0\), for every \(\lambda\le1+\varepsilon_0\),
	\[
	I_\lambda(e)
	=
	\lambda\Psi(e)-\Phi(e)
	\le
	(1+\varepsilon_0)\Psi(e)-\Phi(e)
	=
	I_{1+\varepsilon_0}(e)<0.
	\]
Thus \(e\) is a common endpoint throughout the parameter interval.
\end{remark}
We also have the following result.
\begin{lemma}\label{lem:SMP-and-divergence} 
The family \(I_\lambda=\lambda\Psi-\Phi\) satisfies (SMP).
\end{lemma} 
\begin{proof} 
Lemma \ref{lem:even-condition} shows that \(\Psi\) and \(\Phi\) are even.
Lemma \ref{MP-Q} gives
	\[ 
	\inf_{\|u\|_X=\rho_0} I_{1-\varepsilon_0}(u) \ge\alpha_0>0. 
	\] 
For each \(k\in\mathbb N\), it remains to construct a continuous odd map
	\(\pi_{0,k}\in C(\partial D^k,X),\)
such that
	\begin{equation}\label{smptask}
		\min_{\zeta\in\partial D^k}
		\|\pi_{0,k}(\zeta)\|_X>\rho_0,
		\qquad
		\sup_{\zeta\in\partial D^k}
		I_{1+\varepsilon_0}
		\bigl(\pi_{0,k}(\zeta)\bigr)<0.
	\end{equation}
	
Fix \(k\) and choose linearly independent functions
	\[ 
	\varphi_1,\ldots,\varphi_k \in C_c^\infty(\mathbb R^N). 
	\] 
For \(\zeta=(\zeta_1,\ldots,\zeta_k)\in\partial D^k\), define
	\[ 
	w_\zeta = a_k\sum_{j=1}^{k}\zeta_j\varphi_j, 
	\] 
where \(a_k>0\) is sufficiently small that
	\(\sup_{\zeta\in\partial D^k} \|Dw_\zeta\|_\infty<1.\)
Therefore, for every \(\zeta\in\partial D^k\),
	\(w_\zeta\in D(\Psi)\).
	
The set \(\partial D^k\) is compact. Linear independence of \(\varphi_1,\ldots,\varphi_k\) implies that the map \(\zeta
	\longmapsto
	w_\zeta\)
does not vanish on \(\partial D^k\).
Hence there exists \(M_k>0\) such that
	\[
	\inf_{\zeta\in\partial D^k}
	\int_{\mathbb R^N}|w_\zeta|^p\,dx
	\ge M_k.
	\]
Similarly, if \(\|Dw_\zeta\|_2^2=0,\)
then \(Dw_\zeta=0\) almost everywhere in \(\mathbb R^N\).
Since \(w_\zeta\in C_c^\infty(\mathbb R^N)\), this implies
\(w_\zeta\equiv0\), a contradiction. Therefore,
there exists \(d_k>0\) such that
	\[
	\inf_{\zeta\in\partial D^k}
	\|Dw_\zeta\|_2^2
	\ge d_k.
	\]
On the other hand, \(w_\zeta\in C_c^\infty(\mathbb R^N)\) and
	\(V\in C^1(\mathbb R^N)\), 
Thus there exists
\(A_k,B_k>0\) such that
	\[
	\sup_{\zeta\in\partial D^k}
	\int_{\mathbb R^N}G(Dw_\zeta)\,dx
	\le A_k,
	\qquad
	\sup_{\zeta\in\partial D^k}
	\int_{\mathbb R^N}V(x)w_\zeta^2\,dx
	\le B_k.
	\]
	
For \(T\ge1\), define
	\[ 
	\pi_{k,T_k}(\zeta)(x) = T_kw_\zeta\left(\frac{x}{T_k}\right). 
	\] 
Then \(\pi_{k,T_k}(\zeta)\in D(\Psi)\), and \(\pi_{k,T_k}\) is continuous and odd.
Using \eqref{vtytcd} gives
	\[
	I_{1+\varepsilon_0} \bigl(\pi_{k,T_k}(\zeta)\bigr) 
	\le (1+\varepsilon_0)A_kT_k^N
	+ \frac{1+\varepsilon_0}{2} B_kT_k^{N+2+C_D}
	- \frac{M_k}{p}T_k^{N+p}. 
	\] 
Since \(C_D<p-2\), choose \(T_k\) sufficiently large that
	\[ 
	\sup_{\zeta\in\partial D^k} I_{1+\varepsilon_0} \bigl(\pi_{k,T_k}(\zeta)\bigr)<0. 
	\] 
	
Moreover,
	\[
	\|\pi_{k,T_k}(\zeta)\|_X^2
	\ge
	\|D\pi_{k,T_k}(\zeta)\|_2^2
	=
	T_k^N\|Dw_\zeta\|_2^2
	\ge d_kT_k^N.
	\]
allows us to increase \(T_k\) so that
	\[
	\min_{\zeta\in\partial D^k}
	\|\pi_{k,T_k}(\zeta)\|_X>\rho_0.
	\]
	
Set \(\pi_{0,k}:=\pi_{k,T_k}.\)
Then \(\pi_{0,k}\in C(\partial D^k,X)\) is continuous and odd and satisfies \eqref{smptask}.
This proves (SMP).
\end{proof}

We now establish the properties of \(c_k(\lambda)\) and \(k\in\mathbb N_0=\{0,1,2,\ldots\}\) used to select the parameter sequence \(\lambda\).
\begin{lemma}\label{lem:minimax-values-basic}
Suppose that the family \(I_\lambda=\lambda\Psi-\Phi\)
satisfies the mountain pass conditions in Lemmas \ref{MP-Q} and \ref{lem:SMP-and-divergence}.
For every \(k\in\mathbb N_0\), the minimax values
\(c_k(\lambda)\) satisfy
		\[
		0<\alpha_0\le c_k(\lambda)<+\infty,
		\qquad
		\lambda\in[1-\varepsilon_0,1+\varepsilon_0].
		\]
Moreover, the function \(\lambda\mapsto c_k(\lambda)\)
are nondecreasing on \([1-\varepsilon_0,1+\varepsilon_0]\),
and are therefore differentiable almost everywhere on this interval.
 \end{lemma}

\begin{proof}
By Lemma \ref{MP-Q},
Since \(\|e\|_X>\rho_0\), every \(\gamma\in\Gamma_0\) crosses the sphere
	\[
	\{u\in X:\|u\|_X=\rho_0\}.
	\]
For \(\lambda\ge1-\varepsilon_0\),
	\(I_\lambda(u)\ge I_{1-\varepsilon_0}(u).\)
Therefore,
	\(c_0(\lambda)\ge\alpha_0>0.\)
Conversely, \(t\mapsto te\) belongs to \(\Gamma_0\), so
	\(c_0(\lambda)<+\infty.\)
Consequently,
	\[
	0<\alpha_0\le c_0(\lambda)<+\infty,
	\qquad
	\lambda\in[1-\varepsilon_0,1+\varepsilon_0].
	\]
	
For the symmetric minimax levels, fix \(k\in\mathbb N\).
First, \(\Gamma_k\ne\varnothing\). Define \(\gamma_k(0)=0,\)
and, for \(\zeta\ne0\), put
	\[
	\gamma_k(\zeta)
	=
	|\zeta|\,
	\pi_{0,k}
	\left(\frac{\zeta}{|\zeta|}\right),
	\qquad
	\zeta\in D^k.
	\]
Compactness of \(\pi_{0,k}(\partial D^k)\) gives \(\gamma_k\in C(D^k,X).\)
Moreover, \(\gamma_k\) is odd and \(\gamma_k|_{\partial D^k}=\pi_{0,k}.\)
Thus \(\gamma_k\in\Gamma_k,\)
and hence \(\Gamma_k\ne\varnothing\).
	
Let \(\gamma\in\Gamma_k\). Since \(\gamma\) is odd, \(\gamma(0)=0.\)
Together with
	\[
	\min_{\zeta\in\partial D^k}
	\|\gamma(\zeta)\|_X
	=
	\min_{\zeta\in\partial D^k}
	\|\pi_{0,k}(\zeta)\|_X
	>
	\rho_0,
	\]
this implies that, for each fixed \(\zeta_*\in\partial D^k\), \(\gamma(t\zeta_*)\)
crosses the sphere
	\[
	\{u\in X:\|u\|_X=\rho_0\}.
	\]
There is therefore \(t_\gamma\in(0,1)\) such that \(\|\gamma(t_\gamma\zeta_*)\|_X=\rho_0.\)
For \(\lambda\ge1-\varepsilon_0\),
	\[
	I_\lambda(\gamma(t_\gamma\zeta_*))
	\ge
	I_{1-\varepsilon_0}(\gamma(t_\gamma\zeta_*))
	\ge
	\alpha_0.
	\]
Hence \(\sup_{\zeta\in D^k}I_\lambda(\gamma(\zeta))
	\ge
	\alpha_0.\)
Taking the infimum over \(\gamma\in\Gamma_k\) gives
	\(c_k(\lambda)\ge\alpha_0.\)
	
Conversely, the previously constructed \(\gamma_k\in\Gamma_k\) gives
	\[
	c_k(\lambda)
	\le
	\sup_{\zeta\in D^k}
	I_\lambda(\gamma_k(\zeta)).
	\]
Since \(\gamma_k(D^k)\) is compact in \(X\) and
\(\gamma_k(D^k)\subset D(\Psi)\),
	\(\sup_{\zeta\in D^k}
	I_\lambda(\gamma_k(\zeta))<+\infty.\)
we obtain \(c_k(\lambda)<+\infty.\)
Thus, for every \(k\in\mathbb N\), \(0<\alpha_0\le c_k(\lambda)<+\infty.\)
	
Finally, if \(1-\varepsilon_0\le \lambda_1<\lambda_2\le1+\varepsilon_0,\)
then \(\Psi\ge0\) implies, for every \(u\in X\),
	\[
	I_{\lambda_1}(u)
	=
	\lambda_1\Psi(u)-\Phi(u)
	\le
	\lambda_2\Psi(u)-\Phi(u)
	=
	I_{\lambda_2}(u).
	\]
Since neither \(\Gamma_0\) nor \(\Gamma_k\) depends on \(\lambda\),
taking the minimax values yields, for every \(k\in\mathbb N_0\),
	\[
	c_k(\lambda_1)\le c_k(\lambda_2).
	\]
Consequently the functions \(\lambda\mapsto c_k(\lambda)\)
are nondecreasing on \([1-\varepsilon_0,1+\varepsilon_0]\).
The standard differentiability theorem for monotone functions shows that \(c_k\) are differentiable almost everywhere on the interval.
\end{proof}
\begin{remark}
In applying Lemma \ref{Psexist} to obtain PS sequences, we work at parameters where \(c_k'(\lambda)\) exists.
For the symmetric levels \(c_k(\lambda)\) with \(k\ge1\), we do not construct a separate parameter-dependent critical-point sequence.
After verifying (E), (SMP), \((\Phi_2)\), and (IB),
we apply Lemma \ref{th1.7} directly to obtain infinitely many critical points of \(I_1\).
\end{remark}

We now consider only the critical points associated with \(c_0(\lambda)\).
Choose \(\{\lambda_n\}_{n=1}^\infty\subset(1,1+\epsilon_0]\) such that
\begin{equation}\label{lambdanchoose}
	1<\lambda_n<\min\left\{1+\frac1n,\lambda_{n-1}\right\},
\end{equation}
and \(c_0'(\lambda_n)\) exists. Then \(\lambda_n\downarrow1\).
Fix this sequence \(\{\lambda_n\}\). By Lemma \ref{Psexist}, for every \(n\) there exists a PS sequence for \(I_{\lambda_n}\)
\(\{u_{n,j}\}_{j=1}^{\infty}\subset D(\Psi)\) bounded in \(X\).
Thus, as \(j\to\infty\),
\[
I_{\lambda_n}(u_{n,j})\to c_0(\lambda_n),
\]
and there exists \(\varepsilon_{n,j}\to0\) such that, for every \(v\in D(\Psi)\),
\begin{equation}\label{eq:lamdanps}
	\lambda_n\bigl(\Psi(v)-\Psi(u_{n,j})\bigr)
	-
	\Phi'(u_{n,j})(v-u_{n,j})
	\ge
	-\varepsilon_{n,j}\|v-u_{n,j}\|_X.
\end{equation}

\begin{lemma}\label{lem:low-threshold-compactness}
Fix \(n\in\mathbb N\). Let \(\{u_{n,j}\}_{j=1}^{\infty}\subset D(\Psi)\) be the bounded PS sequence for \(I_{\lambda_n}\) obtained from Lemma \ref{Psexist}, satisfying \(I_{\lambda_n}(u_{n,j})\to c_0(\lambda_n)\)
and \eqref{eq:lamdanps}. Then \(\{u_{n,j}\}_{j=1}^{\infty}\) has a strongly convergent subsequence in \(X\).
\end{lemma}

\begin{proof}
Since \(\{u_{n,j}\}_{j=1}^{\infty}\) is bounded in the reflexive Banach space \(X\), with \(X\) reflexive, there exists \(u_n\in X\) such that, along a subsequence, \(u_{n,j}\rightharpoonup u_n\) in \(X\).
By Lemma \ref{lem:Dpsi-weakly-closed}, \(D(\Psi)\) is weakly closed in \(X\), so \(u_n\in D(\Psi).\)
	
The compact embedding \(X\hookrightarrow L^p(\mathbb R^N),\) from Lemma \ref{space-properties}
gives \(u_{n,j}\to u_n\) in \(L^p(\mathbb R^N)\).
Write
	\(z_{n,j}=u_{n,j}-u_n\), 
Then, in \(L^p(\mathbb R^N)\),
	\(z_{n,j}\to0\). 
	
We improve this strong \(L^p\) convergence to strong convergence in \(X\). Taking \(v=u_n\) in \eqref{eq:lamdanps} gives
	\begin{equation}\label{gdps}
		\lambda_n\bigl(\Psi(u_n)-\Psi(u_{n,j})\bigr)
		-
		\Phi'(u_{n,j})(u_n-u_{n,j})
		\ge
		-\varepsilon_{n,j}\|u_n-u_{n,j}\|_X.
	\end{equation}
Boundedness of \(\{u_{n,j}\}\) in \(X\) and the continuous embedding \(X\hookrightarrow L^p(\mathbb R^N)\) give boundedness of \(\{u_{n,j}\}\) in \(L^p(\mathbb R^N)\). H{\"o}lder's inequality yields
	\[
	\begin{aligned}
		\left|
		\Phi'(u_{n,j})(u_n-u_{n,j})
		\right|
		&=
		\left|
		\int_{\mathbb R^N}
		|u_{n,j}|^{p-2}u_{n,j}(u_n-u_{n,j})\,dx
		\right|\\
		&\le
		\|u_{n,j}\|_p^{p-1}
		\|u_n-u_{n,j}\|_p
		\to0.
	\end{aligned}
	\]
Boundedness of \(\{u_{n,j}\}\) in \(X\) implies
	\(\varepsilon_{n,j}\|u_n-u_{n,j}\|_X\to0.\)
Combining this with \eqref{gdps} gives
	\[
	\limsup_{j\to\infty}\Psi(u_{n,j})\le \Psi(u_n).
	\]
Weak lower semicontinuity of \(\Psi\) now gives
	\(\Psi(u_{n,j})\to\Psi(u_n).\)
Since \(u_{n,j}\rightharpoonup u_n\) in \(X\),
condition \((\Psi_3)\) yields a strongly convergent subsequence of \(\{u_{n,j}\}_{j=1}^{\infty}\) in \(X\), with \(u_{n,j}\to u_n\) in \(X\).
\end{proof}

\begin{theorem}\label{coro:critical-point-lambdan}
For each \(n\in\mathbb N\), the strong limit \(u_n\in D(\Psi)\) is a critical point of \(I_{\lambda_n}\) at \(c_0(\lambda_n)\), namely
	\[
	I_{\lambda_n}(u_n)=c_0(\lambda_n),
	\]
and
	\[
	\lambda_n\bigl(\Psi(v)-\Psi(u_n)\bigr)
	-
	\Phi'(u_n)(v-u_n)
	\ge0,
	\qquad
	\forall v\in D(\Psi).
	\]
\end{theorem}

\begin{proof}
Lemma \ref{lem:low-threshold-compactness} gives \(u_{n,j}\to u_n\in D(\Psi)\) in \(X\).
We show that \(u_n\) is critical for \(I_{\lambda_n}\).
For every \(v\in D(\Psi)\),
	\[
	\Psi(u_{n,j})\to\Psi(u_n),
	\qquad
	\Phi'(u_{n,j})(v-u_{n,j})
	\to
	\Phi'(u_n)(v-u_n). 
	\]
	
Letting \(j\to\infty\) and using \eqref{eq:lamdanps} gives
	\begin{equation}\label{unlim}
		\lambda_n\bigl(\Psi(v)-\Psi(u_n)\bigr)
		-
		\Phi'(u_n)(v-u_n)
		\ge0,
		\qquad
		\forall v\in D(\Psi).
	\end{equation}
Hence \(u_n\) is critical for \(I_{\lambda_n}\).
	
The convergences
	\(I_{\lambda_n}(u_{n,j})\to I_{\lambda_n}(u_n)\)
and \(I_{\lambda_n}(u_{n,j})\to c_0(\lambda_n)\)
imply
	\begin{equation}\label{IC}
		I_{\lambda_n}(u_n)=c_0(\lambda_n).
	\end{equation}
\end{proof}

\section{Regularity of critical points}\label{chap:critical-regularity}

We identify nonsmooth critical points with weak solutions and establish their regularity.

Section 3 gives critical points of \(I_{\lambda_n}\) along a sequence \(\lambda_n\downarrow1\).
Fix \(\lambda>0\) and \(\widehat u\in D(\Psi)\), and suppose that \(\widehat u\) is a nonsmooth critical point of \(I_\lambda\), namely
\[
\lambda\bigl(\Psi(v)-\Psi(\widehat u)\bigr)
-
\Phi'(\widehat u)(v-\widehat u)
\ge0,
\qquad
\forall v\in D(\Psi).
\]
With \eqref{rholu}, the corresponding equation becomes
\begin{equation}
	\operatorname{div}\left(
	\frac{D\widehat u}{\sqrt{1-|D\widehat u|^2}}
	\right)
	+
	\rho_{\lambda,\widehat u}
	=
	0.
\end{equation}

We first relate critical points, local minimizers, and weak solutions, and then obtain \(W^{2,q}_{loc}\) regularity and strict spacelikeness.

\subsection{Local minimizers and weak solutions}
We define local minimizers and show that critical points induce them. We then recall a Born--Infeld regularity and weak-solvability result and check its hypotheses for the source of the present non-autonomous equation.
\begin{definition}\label{def:local-admissible-class}
Let \(\Omega\Subset\mathbb R^N\) be a bounded domain and \(\varphi\in C(\partial\Omega)\).
Define the weakly spacelike class
	\[
	\mathcal Y_{\varphi}(\Omega)
	=
	\left\{
	\psi\in W^{1,\infty}(\Omega) \  \middle|\
|D\psi|\le1\ \text{almost everywhere in }\Omega,
\psi=\varphi\ \text{on }\partial\Omega
	\right\}.
	\]
	
For \(\rho\in L^\infty_{\mathrm{loc}}(\mathbb R^N)\), define
	\begin{equation}
		\mathcal J_{\rho}^{\Omega}(\psi)
		=
		\int_{\Omega}G(D\psi)\,dx
		-
		\int_{\Omega}\rho\psi\,dx,
		\qquad
		\psi\in\mathcal Y_{\varphi}(\Omega).
	\end{equation}
If \(w\in\mathcal Y_{\varphi}(\Omega)\) and
	\begin{equation}
		\mathcal J_{\rho}^{\Omega}(w)
		\le
		\mathcal J_{\rho}^{\Omega}(\psi),
		\qquad
		\forall \psi\in\mathcal Y_{\varphi}(\Omega),
	\end{equation}
then \(w\) minimizes, on \(\Omega\) and with boundary values \(\varphi\), the action
	\(\mathcal J_{\rho}^{\Omega}\).
	
More generally, let \(w\in \operatorname{Lip}(\mathbb R^N)\) and suppose that \(|Dw|\le1\) almost everywhere in \(\mathbb R^N\).
If, for every bounded domain \(\Omega\Subset\mathbb R^N\), the restriction
\(w|_{\Omega}\) minimizes, with boundary values \(w|_{\partial\Omega}\),
\(\mathcal J_{\rho}^{\Omega}\), that is,
	\begin{equation}
		\mathcal J_{\rho}^{\Omega}(w)
		\le
		\mathcal J_{\rho}^{\Omega}(\psi),
		\qquad
		\forall \psi\in\mathcal Y_{w|_{\partial\Omega}}(\Omega),
	\end{equation}
then \(w\) is a local minimizer for the source \(\rho\).
\end{definition}

\begin{remark}
Boundary values are usually interpreted as traces. For insufficiently regular \(\partial\Omega\), we use the continuous boundary-value convention of \cite{RN58},
under which admissible functions extend continuously to $\partial\Omega$ and attain $\varphi$.
\end{remark}
Nonsmooth critical points induce local minimizers.
\begin{lemma}\label{lem:critical-local-minimizer-new}
Let \(\lambda>0\) and let \(\widehat u\in D(\Psi)\) be a nonsmooth critical point of \(I_\lambda\).
Then \(\widehat u\) is a local minimizer for the source \(\rho_{\lambda,\widehat u}\) and action
\(\mathcal J_{\rho_{\lambda,\widehat u}}^\Omega\).
\end{lemma}

\begin{proof}
For \(v\in D(\Psi)\), set \(v_t=(1-t)\widehat u+tv,
	t\in(0,1].\)
Convexity of \(D(\Psi)\) gives \(v_t\in D(\Psi).\)
Taking \(v_t\) in the critical point inequality yields
	\[
	\lambda\bigl(\Psi(v_t)-\Psi(\widehat u)\bigr)
	-
	\Phi'(\widehat u)(v_t-\widehat u)
	\ge0.
	\]
By \(v_t-\widehat u=t(v-\widehat u),\)
this can be written as
	\[
	\lambda\bigl(\mathcal K(v_t)-\mathcal K(\widehat u)\bigr)
	+
	\frac{\lambda}{2}
	\int_{\mathbb R^N}V(x)(v_t^2-\widehat u^2)\,dx
	-
	t\int_{\mathbb R^N}|\widehat u|^{p-2}\widehat u(v-\widehat u)\,dx
	\ge0.
	\]
Convexity of \(\mathcal K\) gives
	\[
	\mathcal K(v_t)-\mathcal K(\widehat u)
	\le
	t\bigl(\mathcal K(v)-\mathcal K(\widehat u)\bigr).
	\]
A direct calculation also gives
	\[
	\frac12\int_{\mathbb R^N}V(x)(v_t^2-\widehat u^2)\,dx
	=
	t\int_{\mathbb R^N}V(x)\widehat u(v-\widehat u)\,dx
	+
	\frac{t^2}{2}
	\int_{\mathbb R^N}V(x)(v-\widehat u)^2\,dx.
	\]
Consequently,
	\[
	\begin{aligned}
		0
		&\le
		\lambda t\bigl(\mathcal K(v)-\mathcal K(\widehat u)\bigr)
		+
		\lambda t\int_{\mathbb R^N}V(x)\widehat u(v-\widehat u)\,dx  \\
		&\quad
		-
		t\int_{\mathbb R^N}|\widehat u|^{p-2}\widehat u(v-\widehat u)\,dx
		+
		\frac{\lambda t^2}{2}
		\int_{\mathbb R^N}V(x)(v-\widehat u)^2\,dx.
	\end{aligned}
	\]
Dividing successively by \(\lambda\) and \(t\) and letting \(t\downarrow0\) yields
	\[
	\mathcal K(v)-\mathcal K(\widehat u)
	-
	\int_{\mathbb R^N}
	\left(
	\lambda^{-1}|\widehat u|^{p-2}\widehat u
	-
	V(x)\widehat u
	\right)
	(v-\widehat u)\,dx
	\ge0.
	\]
By the definition of \(\rho_{\lambda,\widehat u}\), this is
	\begin{equation}\label{eq:K-global-vi-hat-u}
		\mathcal K(v)-\mathcal K(\widehat u)
		-
		\int_{\mathbb R^N}
		\rho_{\lambda,\widehat u}(v-\widehat u)\,dx
		\ge0,
		\qquad
		\forall v\in D(\Psi).
	\end{equation}
	
Now choose a bounded domain
\(\Omega\Subset\mathbb R^N\) and put
	\(\varphi=\widehat u|_{\partial\Omega}.\)
For any \(\psi\in\mathcal Y_{\varphi}(\Omega).\)
glue \(\psi\) to
\(\widehat u\) by defining
	\begin{equation}
		v(x)=
		\begin{cases}
			\psi(x), & x\in\Omega,\\
			\widehat u(x), & x\in\mathbb R^N\setminus\Omega.
		\end{cases}
	\end{equation}
Then \(v\in D(\Psi)\). Taking this \(v\) in \eqref{eq:K-global-vi-hat-u},
and using \(v=\widehat u\) on \(\mathbb R^N\setminus\Omega\), gives
	\begin{equation}
		\int_{\Omega}
		\bigl[
		G(D\psi)-G(D\widehat u)
		\bigr]\,dx
		-
		\int_{\Omega}
		\rho_{\lambda,\widehat u}(\psi-\widehat u)\,dx
		\ge0.
	\end{equation}
which is equivalent to
	\[
	\mathcal J_{\rho_{\lambda,\widehat u}}^\Omega(\widehat u)
	\le
	\mathcal J_{\rho_{\lambda,\widehat u}}^\Omega(\psi).
	\]
Since \(\Omega\Subset\mathbb R^N\) and
\(\psi\in\mathcal Y_{\widehat u|_{\partial\Omega}}(\Omega)\) are arbitrary,
\(\widehat u\) is a local minimizer corresponding to \(\rho_{\lambda,\widehat u}\).
\end{proof}

Under appropriate hypotheses, local minimizers are equivalent to weak solutions of the Born--Infeld equation and enjoy regularity \cite{Byeon2026}.
	\begin{proposition}\label{prop:BI-local-minimizer-regularity}
Let \(u\in\operatorname{Lip}(\mathbb R^N)\) satisfy
		\[
		|Du|\le1
		\quad
		\text{a.e. in }\mathbb R^N,
		\]
and let \(\rho\in L_{\mathrm{loc}}^\infty(\mathbb R^N)\). Suppose that there are a constant
\(c_1>0\) and a function \(h:\mathbb R^N\to(0,\infty)\) such that
		\[
		h(x)\to+\infty,
		\quad
		 |x| \to \infty; 
		\qquad \lim_{|x| \to \infty} \frac{h(x)}{|x|} = 0,
		\]
and
		\begin{equation}\label{eq:BI-growth-condition}
			|u(x)|
			\le
			c_1+|x|-h(x),
			\qquad
			x\in\mathbb R^N.
		\end{equation}
Then the following statements are equivalent:
		\begin{itemize}
\item \(u\) is a weak solution of
			\begin{equation}
				\operatorname{div}\left(
				\frac{Du}{\sqrt{1-|Du|^2}}
				\right)
				+
				\rho
				=
				0
			\end{equation}
;
			
\item \(u\) is a local minimizer for the action with source \(\rho\).
		\end{itemize}
If either holds, then
		\[
		u\in W_{\mathrm{loc}}^{2,q}(\mathbb R^N),
		\qquad
		\forall q\in[2,\infty),
		\]
and \(|Du|<1\) in \(\mathbb R^N\).
	\end{proposition}
We verify these hypotheses for our model, beginning with the source.
\begin{lemma}\label{lem:rho-u-L2loc}
For every \(\lambda>0\) and \(u\in D(\Psi)\cap X\),
	\begin{equation}\label{rholuem}
		\rho_{\lambda,u}\in L^\infty_{loc}(\mathbb R^N)
		\subset
		L^2_{loc}(\mathbb R^N).
	\end{equation}
\end{lemma}

\begin{proof}
Fix a compact set \(K\Subset\mathbb R^N\). Since \(V\in C^1(\mathbb R^N)\),
	\(\|V\|_{L^\infty(K)}<+\infty .\)
Since \(u\in X\hookrightarrow L^\infty_{loc}(\mathbb R^N)\) and
\(u\in H^1_V(\mathbb R^N)\subset H^1(\mathbb R^N)\),
	\(u\in L^\infty(K)\cap L^2(K).\)
	
we first have
	\[
	\|Vu\|_{L^\infty(K)}
	\le
	\|V\|_{L^\infty(K)}\|u\|_{L^\infty(K)}
	<+\infty .
	\]
Also,
	\[
	\bigl\|\lambda^{-1}|u|^{p-2}u\bigr\|_{L^\infty(K)}
	\le
	\lambda^{-1}\|u\|_{L^\infty(K)}^{p-1}
	<+\infty .
	\]
Therefore,
	\[
	\lambda^{-1}|u|^{p-2}u\in L^\infty(K),
	\qquad
	Vu\in L^\infty(K).
	\]
This gives
	\(\rho_{\lambda,u}\in L^\infty(K).\)
	
On the other hand,
	\[
	\|V u\|_{L^2(K)}
	\le
	\|V\|_{L^\infty(K)}\|u\|_{L^2(K)}
	<+\infty .
	\]
Also,
	\[
	\begin{aligned}
		\bigl\|\lambda^{-1}|u|^{p-2}u\bigr\|_{L^2(K)}^2
		&=
		\lambda^{-2}\int_K |u|^{2(p-1)}\,dx  \\
		&\le
		\lambda^{-2}\|u\|_{L^\infty(K)}^{2(p-2)}
		\int_K |u|^2\,dx
		<+\infty .
	\end{aligned}
	\]
Therefore,
	\[
	\lambda^{-1}|u|^{p-2}u\in L^2(K),
	\qquad
	V u\in L^2(K).
	\]
and hence
	\(\rho_{\lambda,u}\in L^2(K).\)
Since \(K\Subset\mathbb R^N\) is arbitrary,
	\[
	\rho_{\lambda,u}\in L^\infty_{loc}(\mathbb R^N)
	\subset
	L^2_{loc}(\mathbb R^N).
	\]
\end{proof}
We can now verify the assumptions of Proposition \ref{prop:BI-local-minimizer-regularity}.
\begin{lemma}\label{lem:verify-BI-regularity-assumptions}
Let \(\lambda>0\) and \(\widehat u\in D(\Psi)\), and set
	\[
	\rho_{\lambda,\widehat u}
	=
	\lambda^{-1}|\widehat u|^{p-2}\widehat u
	-
	V(x)\widehat u.
	\]
Then \(\widehat u\) satisfies the hypotheses of Proposition \ref{prop:BI-local-minimizer-regularity}.
\end{lemma}

\begin{proof}
From \(\widehat u\in D(\Psi)\), we have \(\widehat u\in X,\) and \(|D\widehat u|\le 1\) almost everywhere in \(\mathbb R^N\).
Using \eqref{morreyembed},
we choose a representative of \(\widehat u\),
\(C^{0,\alpha}_{loc}(\mathbb R^N)\) continuous representative of \(\widehat u\). In particular,
for every \(R>0\), \(\widehat u\in L^\infty(B_R).\)
On the other hand, \(\|D\widehat u\|_{L^\infty(B_R)}\le1\) implies
	\(\widehat u\in W^{1,\infty}(B_R)\). 
	
Consequently, for every \(x,y\in B_R\),
	\begin{equation}
		|\widehat u(x)-\widehat u(y)|
		\le
		\|D\widehat u\|_{L^\infty(B_R)}|x-y|
		\le
		|x-y|.
	\end{equation}
Since \(R>0\) is arbitrary, for each \(x,y\in\mathbb R^N\) choose
\(R>|x|+|y|\). Then \(x,y\in B_R\) and hence
	\[
	|\widehat u(x)-\widehat u(y)|
	\le
	|x-y|.
	\]
Thus \(\widehat u\) satisfies \(\widehat u\in \operatorname{Lip}(\mathbb R^N),\) and \(|D\widehat u|\le1\) almost everywhere in \(\mathbb R^N\).
	
The assertion \(\rho_{\lambda,\widehat u}
\in L_{\mathrm{loc}}^\infty(\mathbb R^N)\) follows directly from Lemma \ref{lem:rho-u-L2loc}.
	
Finally, for the growth condition choose \(h(x)=\sqrt{1+|x|}.\)
Clearly,
	\[
	h(x)>0; 
	\qquad
	h(x)\to+\infty,
	\quad
	|x| \to \infty; 
	\qquad \lim_{|x| \to \infty} \frac{h(x)}{|x|} = 0.
	\]
Lemma \ref{space-properties} gives \(X\hookrightarrow C_0(\mathbb R^N) \subset L^\infty(\mathbb R^N)\). Write
	\(M_{\widehat u}\coloneq\|\widehat u\|_\infty<+\infty.\)
Set
	\(c_1=M_{\widehat u}+1,\)
Then, for every \(x\in\mathbb R^N\),
	\[
	c_1+|x|-h(x)
	=
	M_{\widehat u}+1+|x|-\sqrt{1+|x|}
	\ge
	M_{\widehat u}
	\ge
	|\widehat u(x)|.
	\]
Consequently,
	\[
	|\widehat u(x)|
	\le
	c_1+|x|-h(x),
	\qquad
	x\in\mathbb R^N.
	\]
All the assumptions of Proposition \ref{prop:BI-local-minimizer-regularity} are therefore satisfied.
\end{proof}

\subsection{Regularity of nonsmooth critical points}
We establish the equivalence between critical points, local minimizers, and weak solutions, together with \(W_{\mathrm{loc}}^{2,q}\) regularity and strict spacelikeness. We then identify the distributional mean curvature using the geometric formulas of Section 2.

We first give an equivalent formulation of the critical point inequality \eqref{nonsm-critical}.
\begin{lemma}\label{lem:K-vi-to-critical}
Let \(\lambda>0\) and \(\widehat u\in D(\Psi)\). If, for every \(v\in D(\Psi)\),
	\begin{equation}\label{eq:K-vi-before-critical}
		\mathcal K(v)-\mathcal K(\widehat u)
		-
		\int_{\mathbb R^N}
		\rho_{\lambda,\widehat u}(v-\widehat u)\,dx
		\ge0,
	\end{equation}
then \(\widehat u\) satisfies the nonsmooth critical point inequality for \(I_\lambda=\lambda\Psi-\Phi\), namely
	\begin{equation}\label{eq:recover-nonsmooth-critical}
		\lambda\bigl(\Psi(v)-\Psi(\widehat u)\bigr)
		-
		\Phi'(\widehat u)(v-\widehat u)
		\ge0,
		\qquad
		\forall v\in D(\Psi).
	\end{equation}
\end{lemma}

\begin{proof}
For \(v\in D(\Psi)\), the assumption \(\lambda>0\) and \eqref{eq:K-vi-before-critical}
yield
	\begin{equation}\label{eq:K-vi-lambda-expanded}
		\lambda\bigl(\mathcal K(v)-\mathcal K(\widehat u)\bigr)
		+
		\lambda\int_{\mathbb R^N}V(x)\widehat u(v-\widehat u)\,dx
		-
		\int_{\mathbb R^N}|\widehat u|^{p-2}\widehat u(v-\widehat u)\,dx
		\ge0.
	\end{equation}
On the other hand, \(V(x)\ge V_0>0\) gives
	\begin{equation}\label{vinf}
		\frac{\lambda}{2}
		\int_{\mathbb R^N}V(x)(v-\widehat u)^2\,dx
		\ge0.
	\end{equation}
Using
\[
\frac12(v^2-\widehat u^2)
=
\widehat u(v-\widehat u)
+
\frac12(v-\widehat u)^2,
\]
we add the following nonnegative term associated with \eqref{vinf}
term of \eqref{eq:K-vi-lambda-expanded} to the left-hand side of \eqref{eq:K-vi-lambda-expanded} to obtain
	\[
	\lambda\bigl(\mathcal K(v)-\mathcal K(\widehat u)\bigr)
	+
	\frac{\lambda}{2}
	\int_{\mathbb R^N}V(x)(v^2-\widehat u^2)\,dx
	-
	\int_{\mathbb R^N}|\widehat u|^{p-2}\widehat u(v-\widehat u)\,dx
	\ge0.
	\]
Since
	\[
	\lambda\bigl(\Psi(v)-\Psi(\widehat u)\bigr)
	=
	\lambda\bigl(\mathcal K(v)-\mathcal K(\widehat u)\bigr)
	+
	\frac{\lambda}{2}
	\int_{\mathbb R^N}V(x)(v^2-\widehat u^2)\,dx.
	\]
and
	\[
	\Phi'(\widehat u)(v-\widehat u)
	=
	\int_{\mathbb R^N}
	|\widehat u|^{p-2}\widehat u(v-\widehat u)\,dx.
	\]
Consequently,
	\[
	\lambda\bigl(\Psi(v)-\Psi(\widehat u)\bigr)
	-
	\Phi'(\widehat u)(v-\widehat u)
	\ge0.
	\]
with \(v\in D(\Psi)\) arbitrary, this proves \eqref{eq:recover-nonsmooth-critical}.
\end{proof}

\begin{theorem}\label{thm:critical-minimizer-weak-equivalence}\label{thm:critical-to-weak-regular}
Let \(\lambda>0\) and \(\widehat u\in D(\Psi)\).
Then the following three statements are equivalent:
	\begin{itemize}
\item \(\widehat u\) is a nonsmooth critical point of \(I_\lambda=\lambda\Psi-\Phi\);
		
\item \(\widehat u\) is a local minimizer, for the source \(\rho_{\lambda,\widehat u}\), of
\(\mathcal J_{\rho_{\lambda,\widehat u}}^\Omega\);
		
\item \(\widehat u\) is a weak solution in the sense of Definition \ref{def:ws}
, that is,
		\[
		\operatorname{div}\left(
		\frac{D\widehat u}{\sqrt{1-|D\widehat u|^2}}
		\right)
		-
		V(x)\widehat u
		+
		\lambda^{-1}|\widehat u|^{p-2}\widehat u
		=
		0
		\]
in distributions.
	\end{itemize}
If any of these holds, then, for every \(q\in[2,\infty)\), \(\widehat u\in W_{\mathrm{loc}}^{2,q}(\mathbb R^N),\)
and, in \(\mathbb R^N\),
	\(|D\widehat u|<1.\)
The associated mean curvature also satisfies, in distributions,
	\begin{equation}\label{eq:H-rho-lambda-regularity}
		H_{\lambda,\widehat u}
		=
		V(x)\widehat u
		-
		\lambda^{-1}|\widehat u|^{p-2}\widehat u
		=
		-\rho_{\lambda,\widehat u}.
	\end{equation}
\end{theorem}

\begin{proof}
Suppose first that \(\widehat u\) is critical for \(I_\lambda\). Lemma
\ref{lem:critical-local-minimizer-new} shows that \(\widehat u\) is a local minimizer for the source
\(\rho_{\lambda,\widehat u}\) and action
\(\mathcal J_{\rho_{\lambda,\widehat u}}^\Omega\).
	
By Lemma
\ref{lem:verify-BI-regularity-assumptions}, whenever
\(\widehat u\in D(\Psi)\),
all the hypotheses of Proposition \ref{prop:BI-local-minimizer-regularity} hold. That proposition shows that
\(\widehat u\) is a local minimizer for \(\rho_{\lambda,\widehat u}\) if and only if
\(\widehat u\) weakly solves
	\[
	\operatorname{div}\left(
	\frac{D\widehat u}{\sqrt{1-|D\widehat u|^2}}
	\right)
	+
	\rho_{\lambda,\widehat u}
	=
	0
	\]
.

It remains to recover the critical point inequality from weak solvability or local minimality,
\eqref{nonsm-critical}. Since Proposition
\ref{prop:BI-local-minimizer-regularity} gives their equivalence, we start from local minimality.

Let \(\widehat u\) be a local minimizer for \(\rho_{\lambda,\widehat u}\). By Lemma \ref{lem:K-vi-to-critical}, it suffices to prove
\begin{equation}\label{eq:global-K-vi-from-minimizer}
	\mathcal K(v)-\mathcal K(\widehat u)
	-
	\int_{\mathbb R^N}
	\rho_{\lambda,\widehat u}(v-\widehat u)\,dx
	\ge0,
	\qquad
	\forall v\in D(\Psi).
\end{equation}
.

For \(v\in D(\Psi)\), write \(\widehat w=v-\widehat u.\)
Since \(v,\widehat u\in D(\Psi)\subset X\hookrightarrow C_0(\mathbb R^N),\)
we have
\[
\widehat w\in C_0(\mathbb R^N); 
\qquad
\widehat w(x)\to0, 
\quad
 |x|\to+\infty.
\]
Choose a sequence \(\varepsilon_j\downarrow0\) and define the truncations
\[
T_j(s)
=
\begin{cases}
	s-\varepsilon_j, & s>\varepsilon_j,\\
	0, & |s|\le\varepsilon_j,\\
	s+\varepsilon_j, & s<-\varepsilon_j.
\end{cases}
\]
Set
\begin{equation}\label{eq:vj-truncation}
	v_j
	=
	\widehat u+T_j(v-\widehat u).
\end{equation}
Equivalently, \(v_j
=
\max\{\widehat u,v-\varepsilon_j\}
+
\min\{\widehat u,v+\varepsilon_j\}
-
\widehat u. \)

We show that \(v_j\) is an admissible local comparison. By
\(\widehat w\in C_0(\mathbb R^N)\), the set
\(\{x\in\mathbb R^N: |\widehat w(x)|\ge\varepsilon_j\}\)
is compact, whereas
\[
\operatorname{supp}(v_j-\widehat u)
=
\operatorname{supp}T_j(\widehat w)
\subset
\{x\in\mathbb R^N: |\widehat w(x)|\ge\varepsilon_j\}.
\]
Thus \(v_j-\widehat u\) has compact support. Moreover, for every \(R>0\), \(\widehat w=v-\widehat u\in W^{1,\infty}(B_R)\).
Next, \(T_j\) is Lipschitz and \(T_j(0)=0\). Hence \(T_j(\widehat w)\in W^{1,\infty}(B_R).\). Since \(R>0\) is arbitrary,
\(v_j\in W_{\mathrm{loc}}^{1,\infty}(\mathbb R^N)\), and almost everywhere,
\[
Dv_j
=
D\widehat u
+
T_j'(\widehat w)D\widehat w.
\]
that is,
\begin{equation}\label{dvj}
	Dv_j=
	\begin{cases}
		D\widehat u, & |v-\widehat u|\le\varepsilon_j,\\
		Dv, & |v-\widehat u|>\varepsilon_j.
	\end{cases}
\end{equation}
Thus \(|Dv_j|\le1\) almost everywhere in \(\mathbb R^N\).
We also have \(|v_j|\le |\widehat u|+|v|.\)
We have
\[
\int_{\mathbb R^N}V(x)v_j^2\,dx
\le
2\int_{\mathbb R^N}V(x)\widehat u^2\,dx
+
2\int_{\mathbb R^N}V(x)v^2\,dx<+\infty,
\]
and
\[
\|Dv_j\|_2+\|Dv_j\|_{p_1}<+\infty.
\]
and therefore \(v_j\in X,\)\(\Psi(v_j)<+\infty, \), so
\(v_j\in D(\Psi).\)

Choose a bounded domain \(\Omega_j\Subset\mathbb R^N\) such that
\(\operatorname{supp}(v_j-\widehat u)\Subset\Omega_j.\)
Then \(v_j=\widehat u\) on \(\partial\Omega_j\) and \(v_j|_{\Omega_j}
\in
\mathcal Y_{\widehat u|_{\partial\Omega_j}}(\Omega_j).\)
Since \(\widehat u\) is a local minimizer for \(\rho_{\lambda,\widehat u}\),
\(\mathcal J_{\rho_{\lambda,\widehat u}}^{\Omega_j}(\widehat u)
\le
\mathcal J_{\rho_{\lambda,\widehat u}}^{\Omega_j}(v_j), \)
that is,
\[
\int_{\Omega_j}
\bigl[
G(Dv_j)-G(D\widehat u)
\bigr]\,dx
-
\int_{\Omega_j}
\rho_{\lambda,\widehat u}(v_j-\widehat u)\,dx
\ge0.
\]
and since \(v_j=\widehat u\) on \(\mathbb R^N\setminus\Omega_j\), the inequality can be written over the whole space as
\begin{equation}\label{eq:global-K-vi-vj}
	\mathcal K(v_j)-\mathcal K(\widehat u)
	-
	\int_{\mathbb R^N}
	\rho_{\lambda,\widehat u}(v_j-\widehat u)\,dx
	\ge0.
\end{equation}

Let \(j\to\infty\). By construction, for every \(x\in\mathbb R^N\), \(v_j(x)\to v(x)\) holds in \(\mathbb R^N\).
Also \eqref{dvj} implies, almost everywhere, \(Dv_j(x)\to Dv(x).\)
Since \(Dv_j\) equals either \(D\widehat u\) or \(Dv\) almost everywhere,
\eqref{foundation} gives, in \(\mathbb R^N\),
\[
\begin{aligned}
	\left|
	G(Dv_j)-G(D\widehat u)
	\right|
	&=
	\left|
	\sqrt{1-|Dv_j|^2}
	-
	\sqrt{1-|D\widehat u|^2}
	\right|  \\
	&\le
	1-\sqrt{1-|Dv|^2}
	+
	1-\sqrt{1-|D\widehat u|^2}  \\
	&=
	G(Dv)+G(D\widehat u)\\
	&\le
	|Dv|^2+|D\widehat u|^2.
\end{aligned}
\]
almost everywhere.
Dominated convergence yields
\begin{equation}
	\mathcal K(v_j)\to \mathcal K(v).
\end{equation}

For the source integral in \eqref{eq:global-K-vi-vj}, first observe that
\(\rho_{\lambda,\widehat u}(v-\widehat u)\in L^1(\mathbb R^N).\)
From
\[
\rho_{\lambda,\widehat u}
=
\lambda^{-1}|\widehat u|^{p-2}\widehat u
-
V(x)\widehat u, 
\]
we estimate the two terms separately. For the first, \(X\hookrightarrow L^p(\mathbb R^N)\) gives
\[
\int_{\mathbb R^N}
|\widehat u|^{p-1}|v-\widehat u|\,dx
\le
\|\widehat u\|_p^{p-1}
\|v-\widehat u\|_p
<+\infty.
\]
For the second, \(v,\widehat u\in X\) yields
\[
\int_{\mathbb R^N}V(x)(v-\widehat u)^2\,dx
\le
2\int_{\mathbb R^N}V(x)v^2\,dx
+
2\int_{\mathbb R^N}V(x)\widehat u^2\,dx
<+\infty.
\]
and Cauchy--Schwarz gives
\[
\int_{\mathbb R^N}
V(x)|\widehat u||v-\widehat u|\,dx
\le
\left(
\int_{\mathbb R^N}V(x)\widehat u^2\,dx
\right)^{1/2}
\left(
\int_{\mathbb R^N}V(x)(v-\widehat u)^2\,dx
\right)^{1/2}
<+\infty.
\]
Consequently,
\(\rho_{\lambda,\widehat u}(v-\widehat u)\in L^1(\mathbb R^N).\)

On the other hand, \(|v_j-\widehat u|
=
|T_j(v-\widehat u)|
\le
|v-\widehat u|.\)
Therefore,
\[
|\rho_{\lambda,\widehat u}(v_j-\widehat u)|
\le
|\rho_{\lambda,\widehat u}|\,|v-\widehat u|
\in L^1(\mathbb R^N).
\]
so dominated convergence gives
\[
\int_{\mathbb R^N}
\rho_{\lambda,\widehat u}(v_j-\widehat u)\,dx
\to
\int_{\mathbb R^N}
\rho_{\lambda,\widehat u}(v-\widehat u)\,dx.
\]
Thus, in \eqref{eq:global-K-vi-vj},
letting \(j\to\infty\) yields
\[
\mathcal K(v)-\mathcal K(\widehat u)
-
\int_{\mathbb R^N}
\rho_{\lambda,\widehat u}(v-\widehat u)\,dx
\ge0.
\]
This proves \eqref{eq:global-K-vi-from-minimizer}.

We therefore have
\begin{equation}
	\lambda\bigl(\Psi(v)-\Psi(\widehat u)\bigr)
	-
	\Phi'(\widehat u)(v-\widehat u)
	\ge0, 
	\qquad
	\forall v\in D(\Psi)
\end{equation}
and \(\widehat u\) is a nonsmooth critical point of
\(I_\lambda\).
The equivalence of weak solutions, local minimizers, and critical points is established.

Finally, Proposition
\ref{prop:BI-local-minimizer-regularity} and the preceding proof identify the critical point
\(\widehat u\) as the corresponding local minimizer and weak solution. Consequently,
\begin{equation}\label{ucritical}
	\widehat u\in W_{\mathrm{loc}}^{2,q}(\mathbb R^N),
	\qquad
	\forall q\in[2,\infty),
\end{equation}
and, in \(\mathbb R^N\),
\begin{equation}\label{ducr}
	|D\widehat u|<1.
\end{equation}

In Section 2,
\eqref{H2} was derived as a geometric identity in the smooth setting \(u\in C^2(\mathbb R^N)\) with \(|Du|<1\).
\[
H
=
\operatorname{div}\left(
\frac{Du}{\sqrt{1-|Du|^2}}
\right).
\]
For the critical point \(\widehat u\), Proposition
\ref{prop:BI-local-minimizer-regularity} establishes
\eqref{ucritical} and \eqref{ducr}.
We therefore define, in distributions,
\[
H_{\lambda,\widehat u}
\coloneq
\operatorname{div}\left(
\frac{D\widehat u}{\sqrt{1-|D\widehat u|^2}}
\right).
\]
We check that this is well defined and satisfies
\[
H_{\lambda,\widehat u}
=
-\rho_{\lambda,\widehat u}
=
V(x)\widehat u
-
\lambda^{-1}|\widehat u|^{p-2}\widehat u
\]
in distributions.

Fix a compact set \(K\Subset\mathbb R^N\) and \(q>N\).
The inclusion \(\widehat u\in W_{\mathrm{loc}}^{2,q}(\mathbb R^N)\) and Morrey's embedding give
\begin{equation}
	\widehat u\in C_{\mathrm{loc}}^{1,\alpha}(\mathbb R^N),
	\qquad
	\alpha=1-\frac{N}{q}.
\end{equation}
Hence \(D\widehat u\) is continuous on \(K\). By \eqref{ducr},
there exists \(\delta_K>0\) such that \(|D\widehat u|\le1-\delta_K\) on \(K\).
Therefore,
\[
\frac{D\widehat u}{\sqrt{1-|D\widehat u|^2}}
\in L^\infty(K;\mathbb R^N)
\subset L^1(K;\mathbb R^N).
\]
Since \(K\Subset\mathbb R^N\) is arbitrary,
\[
\frac{D\widehat u}{\sqrt{1-|D\widehat u|^2}}
\in L_{\mathrm{loc}}^1(\mathbb R^N;\mathbb R^N).
\]
Thus, for every
\(\eta\in C_c^\infty(\mathbb R^N)\), we have
\begin{equation}\label{weak-h}
	\left\langle H_{\lambda,\widehat u},\eta\right\rangle
	=
	\left\langle
	\operatorname{div}\left(
	\frac{D\widehat u}{\sqrt{1-|D\widehat u|^2}}
	\right),
	\eta
	\right\rangle
	=
	-
	\int_{\mathbb R^N}
	\frac{D\widehat u\cdot D\eta}
	{\sqrt{1-|D\widehat u|^2}}\,dx.
\end{equation}

The weak formulation gives, for every
\(\eta\in C_c^\infty(\mathbb R^N)\),
\[
\int_{\mathbb R^N}
\frac{D\widehat u\cdot D\eta}
{\sqrt{1-|D\widehat u|^2}}\,dx
=
\int_{\mathbb R^N}
\rho_{\lambda,\widehat u}\eta\,dx.
\]
Substituting \eqref{weak-h} yields
\[
\left\langle H_{\lambda,\widehat u},\eta\right\rangle
=
-
\int_{\mathbb R^N}
\rho_{\lambda,\widehat u}\eta\,dx.
\]
As \(\eta\in C_c^\infty(\mathbb R^N)\) is arbitrary, the expression \(\rho_{\lambda,\widehat u}\) gives
\begin{equation}\label{eq:H-rho-lambda-proof}
	H_{\lambda,\widehat u}
	=
	-\rho_{\lambda,\widehat u}
	=
	V(x)\widehat u
	-
	\lambda^{-1}|\widehat u|^{p-2}\widehat u. 
\end{equation}
in distributions.
\end{proof}

\begin{remark}\label{rem:positive-truncation-regularity}
The argument does not depend on the precise nonlinearity on the negative half-line.
If we define
	\[
	\Phi^+(u)
	=
	\frac1p\int_{\mathbb R^N}(u^+)^p\,dx,
	\qquad
	I_\lambda^+(u)
	=
	\lambda\Psi(u)-\Phi^+(u),
	\]
then every nonsmooth critical point \(\tilde{u}\in D(\Psi)\) of \(I_\lambda^+\)
weakly solves the truncated equation
	\[
	\operatorname{div}
	\left(
	\frac{Du}{\sqrt{1-|Du|^2}}
	\right)
	-
	V(x)u
	+
	\lambda^{-1}(u^+)^{p-1}
	=
	0
	\]
Moreover, for every \(q>N\),
	\(\tilde{u}\in W_{\mathrm{loc}}^{2,q}(\mathbb R^N),\)
and \(|D\tilde{u}|<1\) in \(\mathbb R^N\).
\end{remark}

\section{Proofs of the main theorems}
We have verified the functional assumptions, mountain pass geometry, and regularity of critical points. In particular, Section 3 shows that
\(I_\lambda=\lambda\Psi-\Phi\) satisfy
\((\Psi_1)\), \((\Psi_2)\), \((\Psi_3)\), \((\Phi_1)\), \((\Phi_2)\), 
and that evenness \((E)\) and uniform symmetric mountain pass geometry \((SMP)\) hold.
Section 4 identifies every critical point of \(I_\lambda\) with a weak solution of the parameter-dependent equation,
with local \(W^{2,q}\) regularity and strict spacelikeness.

We now pass from the parameter-dependent problem to the original equation. An energy identity establishes (IB), completing the hypotheses of the multiplicity theorem. The truncated functional \(I_\lambda^+=\lambda\Psi-\Phi^+\)
then produces a nontrivial positive weak solution. The multiplicity result
\(\ref{th1.7}\) yields infinitely many nontrivial critical points of \(I_1\), hence weak solutions by Section 4. Finally, local Schauder estimates give classical regularity.

\subsection{Uniform bounds for parameter-dependent critical points}
We first establish an identity for every critical point. It will be used to verify (IB).
\begin{lemma}\label{lem:nonautonomous-pohozaev}
Let \(\lambda>0\) and let \(u\in D(\Psi)\) be a critical point of \(I_\lambda\).
Write
	\[
	K(u)=\int_{\mathbb R^N}G(Du)\,dx,
	\qquad
	B(u)=\int_{\mathbb R^N}V(x)u^2\,dx,
	\]
	\[
	C(u)=\int_{\mathbb R^N}|u|^p\,dx,
	\qquad
	D(u)=\int_{\mathbb R^N}(x\cdot\nabla V(x))u^2\,dx.
	\]
Then
	\begin{equation}\label{eq:nonauto-pohozaev}
	\lambda
	\left[
	NK(u)+\frac{N+2}{2}B(u)+\frac12D(u)
	\right]
	=
	\frac{N+p}{p}C(u).
    \end{equation}
\end{lemma}

\begin{proof}
We use compactly supported tests to justify the identity without differentiating a scaling path in the duality pairing.
By Theorem~\ref{thm:critical-minimizer-weak-equivalence}, the critical point is a strictly spacelike weak solution and belongs to $W^{2,q}_{\rm loc}$ for every finite $q\ge2$.
Morrey's embedding and local Schauder estimates therefore give $u\in C^{2,\alpha}_{\rm loc}$ for some $\alpha\in(0,1)$, by the local argument detailed in Theorem~\ref{thm:main-classical-regularity} below.
This regularity argument uses only the weak equation and strict spacelikeness, not the present identity.
Set
\[
A=\frac{Du}{\sqrt{1-|Du|^2}},\qquad T=\int_{\mathbb R^N}A\cdot Du\,dx.
\]
First we prove that $T$ is finite. Taking $v=(1-h)u$, $0<h<1$, in the critical point inequality gives
\[
\lambda\frac{\Psi(u)-\Psi((1-h)u)}h\le C(u).
\]
The kinetic difference quotients are nonnegative and converge pointwise to $A\cdot Du$. The potential difference quotient tends to $B(u)$.
Fatou's lemma yields
\begin{equation}\label{eq:arxiv-flux-integrability}
\lambda\bigl(T+B(u)\bigr)\le C(u)<\infty.
\end{equation}
Moreover, $u\in X\hookrightarrow C_0(\mathbb R^N)$ and splitting the domain into $\{|Du|\le1/2\}$ and its complement gives
\[
|u||A|\le2|u||Du|+2\|u\|_\infty A\cdot Du\in L^1(\mathbb R^N).
\]
Choose $\chi_R\in C_c^\infty(\mathbb R^N)$ equal to one on $B_R$, vanishing outside $B_{2R}$, with $0\le\chi_R\le1$ and $|D\chi_R|\le c/R$.
Testing the weak equation with $\chi_Ru$ and letting $R\to\infty$ gives
\begin{equation}\label{eq:arxiv-amplitude-identity}
\lambda\bigl(T+B(u)\bigr)=C(u).
\end{equation}
Indeed, the cutoff error is bounded by $(c/R)\int_{B_{2R}\setminus B_R}|u||A|$, which tends to zero.

Next test with $\chi_R x\cdot Du$, which is a compactly supported Lipschitz function by local classical regularity. The relation
\[
A\cdot D(x\cdot Du)=A\cdot Du+x\cdot D\bigl(G(Du)\bigr)
\]
and integration by parts give, as $R\to\infty$,
\begin{equation}\label{eq:arxiv-cutoff-pohozaev}
\lambda\left[NK(u)+\frac N2B(u)+\frac12D(u)-T\right]=\frac NpC(u).
\end{equation}
All boundary terms tend to zero: their absolute values are bounded by a fixed constant times
\[
\int_{B_{2R}\setminus B_R}\bigl(G(Du)+Vu^2+|u|^p+A\cdot Du\bigr)\,dx.
\]
Here $|x||D\chi_R|$ is bounded on the annulus and
$|(x\cdot\nabla V)u^2|\le C_DVu^2$ is integrable.
For the nonlinear term we used
$|u|^{p-2}u\,x\cdot Du=x\cdot D(|u|^p/p)$.
Finally, eliminating $T$ between \eqref{eq:arxiv-amplitude-identity} and \eqref{eq:arxiv-cutoff-pohozaev} yields
\[
\lambda\left[NK(u)+\frac{N+2}{2}B(u)+\frac12D(u)\right]=\frac{N+p}{p}C(u),
\]
as claimed.
\end{proof}

The identity yields uniform bounds for critical points.
\begin{lemma}\label{lem:IB-condition}
The family \(I_\lambda=\lambda\Psi-\Phi\) satisfies (IB).
\end{lemma}

\begin{proof}
Fix \(0<a\le b<+\infty,\) and \(c\in\mathbb R.\)
If \(\mathcal K_c^{[a,b]}=\varnothing\), there is nothing to prove.
Otherwise, let \(u\in\mathcal K_c^{[a,b]}\). There exists
\(\lambda\in[a,b]\) such that \(u\) is critical for \(I_\lambda\), \(u\in D(\Psi)\)
and \(I_\lambda(u)\le c.\)
By Lemma \ref{lem:nonautonomous-pohozaev}, \(u\) satisfies
	\begin{equation}\label{tem-poho}
		\lambda
		\left[
		NK(u)
		+
		\frac{N+2}{2}B(u)
		+
		\frac12D(u)
		\right]
		=
		\frac{N+p}{p}C(u).
	\end{equation}
On the other hand,
	\[
	I_\lambda(u)
	=
	\lambda K(u)
	+
	\frac{\lambda}{2}B(u)
	-
	\frac1pC(u).
	\]
Eliminating \(C(u)\) by means of \eqref{tem-poho} gives
	\[
	I_\lambda(u)
	=
	\frac{\lambda}{N+p}
	\left[
	pK(u)
	+
	\frac{p-2}{2}B(u)
	-
	\frac12D(u)
	\right].
	\]
Using \eqref{V-derivative-growth}, we obtain
	\[
	D(u)
	=
	\int_{\mathbb R^N}(x\cdot DV(x))u^2\,dx
	\le
	C_D
	\int_{\mathbb R^N}V(x)u^2\,dx
	=
	C_D B(u).
	\]
Consequently,
	\[
	I_\lambda(u)
	\ge
	\frac{\lambda}{N+p}
	\left[
	pK(u)
	+
	\frac{p-2-C_D}{2}B(u)
	\right].
	\]
Since \(\lambda\ge a>0\) and \(p-2-C_D>0,\)
there is a constant \(C_a>0\), depending only on \(a,N,p,C_D\), such that
	\[
	I_\lambda(u)
	\ge
	C_a\bigl(K(u)+B(u)\bigr).
	\]
	
If \(c<0\), no critical point can satisfy
\(I_\lambda(u)\le c\). Thus
	\(\mathcal K_c^{[a,b]}=\varnothing\).
If \(c\ge0\), the inequality \(I_\lambda(u)\le c\) gives
	\[
	K(u)+B(u)\le \frac{c}{C_a}.
	\]
Therefore,
	\[
	\|Du\|_2^2
	\le
	2K(u)
	\le
	\frac{2c}{C_a},
	\]
and
	\[
	\int_{\mathbb R^N}V(x)u^2\,dx
	=
	B(u)
	\le
	\frac{c}{C_a}.
	\]
and hence
	\[
	\|u\|_Y^2
	=
	\|Du\|_2^2
	+
	\int_{\mathbb R^N}V(x)u^2\,dx
	\le
	C.
	\]
Furthermore, \(|Du|\le1\) almost everywhere and \(p_1>2\) imply
	\[
	\|Du\|_{p_1}^{p_1}
	=
	\int_{\mathbb R^N}|Du|^{p_1}\,dx
	\le
	\int_{\mathbb R^N}|Du|^2\,dx
	\le C.
	\]
Hence \(\|Du\|_{p_1}^2\le C.\)
Combining the estimates,
	\[
	\|u\|_X^2
	=
	\|u\|_Y^2+\|Du\|_{p_1}^2
	\le C,
	\]
where \(C\) is independent of \(u\) and \(\lambda\).
Thus \(\mathcal K_c^{[a,b]}\) is bounded in \(X\).
This proves (IB).
\end{proof}
\begin{remark}\label{rem:critical-points-uniform-bound}
Condition (IB) gives uniform bounds for the critical-point sequence from Lemma \ref{Psexist}
\(X\) for the critical-point sequence obtained in Lemma \(X\).
	
Choose \(\lambda_n\) as in \eqref{lambdanchoose} of Section 3.
For \(k=0\), Lemma \ref{Psexist} gives
a bounded PS sequence for \(I_{\lambda_n}\) at level \(c_0(\lambda_n)\),
	\begin{equation*}
		\left\{u_{n,j}^{(0)}\right\}_{j=1}^{\infty}
		\subset D(\Psi).
	\end{equation*}
whose limiting critical point \(u_n^{(0)}\in D(\Psi),\) satisfies
	\[
	I_{\lambda_n}\left(u_n^{(0)}\right)=c_0(\lambda_n).
	\]
Since \(c_0\) is nondecreasing and
	\[
	c_0(\lambda_n)\le c_0(1+\varepsilon_0)<+\infty,
	\]
(IB) implies
	\begin{equation}\label{supunw}
		\sup_{n\in\mathbb N}
		\left\|u_n^{(0)}\right\|_X<+\infty.
	\end{equation}
Similarly, for each fixed \(k\in\mathbb N\), if
\(c_k'(\lambda_n)\) exists, then
	\begin{equation}\label{supunkw}
		\sup_{n\in\mathbb N}
		\left\|u_n^{(k)}\right\|_X<+\infty.
	\end{equation}
\end{remark}

\subsection{The parameter limit and classical solutions}
We obtain a nontrivial critical point of \(I_1\).
\begin{lemma}\label{thm:critical-point-lambda-one}
The functional \(I_1\) has a nontrivial critical point \(u\in D(\Psi)\) satisfying
	\[
	I_1(u)=c_0(1).
	\]
\end{lemma}

\begin{proof}
We first prove strong subsequential convergence of \(\{u_n\}\). By \eqref{supunw},
\(\{u_n\}\) is uniformly bounded in \(X\). Reflexivity of \(X\) gives \(u\in X\) and a subsequence, still denoted by \(\{u_n\}\), such that \(u_n\rightharpoonup u\) in \(X\).
Weak closedness of \(D(\Psi)\) in \(X\) gives \(u\in D(\Psi).\)
To improve this to strong convergence, condition \((\Psi_3)\) reduces the problem to proving \(\Psi(u_n)\to\Psi(u)\).
Lemma \ref{space-properties} gives the compact embedding
	\(X\hookrightarrow L^p(\mathbb R^N),\)
so, after passing to a subsequence of \(u_n\) still denoted by \(u_n\), we have \(u_n\to u\) in \(L^p(\mathbb R^N)\).

Since \(u_n\) is critical for \(I_{\lambda_n}\), for every
\(v\in D(\Psi)\), we have
	\[
	\lambda_n\bigl(\Psi(v)-\Psi(u_n)\bigr)
	-
	\Phi'(u_n)(v-u_n)
	\ge0.
	\]
As \(u\in D(\Psi)\), taking \(v=u\) yields
	\[
	\lambda_n\bigl(\Psi(u)-\Psi(u_n)\bigr)
	\ge
	\Phi'(u_n)(u-u_n).
	\]
H{\"o}lder's
inequality gives
	\[
	\begin{aligned}
		\left|
		\Phi'(u_n)(u-u_n)
		\right|
		&=
		\left|
		\int_{\mathbb R^N}|u_n|^{p-2}u_n(u-u_n)\,dx
		\right|  \\
		&\le
		\|u_n\|_p^{p-1}\|u-u_n\|_p
		\to0.
	\end{aligned}
	\]
Since \(\lambda_n\ge1\), it follows that \(\limsup_{n\to\infty}\Psi(u_n)\le \Psi(u).\)
Weak lower semicontinuity of \(\Psi\) in \(X\) gives
	\(\Psi(u)\le \liminf_{n\to\infty}\Psi(u_n).\)
Consequently \(\Psi(u_n)\to\Psi(u).\)
Condition \((\Psi_3)\) now yields \(u_n\to u\) in \(X\).
	
We next show that \(u\) is critical for \(I_1\).
Since \(\lambda_n\to1,\) and \(u_n\to u\) in \(X\),
and since \(\Psi(u_n)\to\Psi(u)\) and \(\Phi'(u_n)(v-u_n)\to
\Phi'(u)(v-u)\), letting \(n\to\infty\) gives
	\[
	\Psi(v)-\Psi(u)
	-
	\Phi'(u)(v-u)
	\ge0,
	\qquad
	\forall v\in D(\Psi).
	\]
Hence \(u\) is critical for \(I_1\).
	
Moreover, \(I_{\lambda_n}(u_n)=c_0(\lambda_n)\)
and \(I_1(u_n)
	=
	I_{\lambda_n}(u_n)-(\lambda_n-1)\Psi(u_n),\)
imply
	\[
	I_1(u_n)
	=
	c_0(\lambda_n)-(\lambda_n-1)\Psi(u_n).
	\]
Boundedness of \(\{\Psi(u_n)\}\) gives \((\lambda_n-1)\Psi(u_n)\to0.\)
By Lemma \ref{lem:c0-right-continuity},
\(c_0\) is right-continuous at \(\lambda=1\), and hence
	\(	c_0(\lambda_n)\to c_0(1).\)
Consequently,
	\[
	I_1(u)
	=
	\lim_{n\to\infty}I_1(u_n)
	=
	\lim_{n\to\infty}
	\left[
	c_0(\lambda_n)-(\lambda_n-1)\Psi(u_n)
	\right]
	=
	c_0(1).
	\]
	
Finally, mountain pass geometry and Lemma \ref{MP-Q} give
	\(c_0(1)\ge \alpha_0>0.\)
we obtain \(I_1(u)=c_0(1)>0=I_1(0).\)
Hence \(u\not\equiv0\)
so the critical point is nontrivial.
\end{proof}
A corresponding argument for the positive truncation gives a positive weak solution.
\begin{theorem}\label{solution-positive}
The original problem \eqref{biff} admits a nontrivial positive weak solution \(w\in Y\).
\end{theorem}
\begin{proof}
Define the truncated functional
	\[
	\Phi^+(u)
	=
	\frac1p\int_{\mathbb R^N}(u^+)^p\,dx,
	\qquad
	I_\lambda^+(u)
	=
	\lambda\Psi(u)-\Phi^+(u),
	\quad
	u^+=\max\{u,0\}.
	\]
By the embedding in Lemma \ref{space-properties},
\(X\hookrightarrow L^p(\mathbb R^N)\), we have
\(\Phi^+\in C^1(X,\mathbb R)\) and
	\[
	(\Phi^+)'(u)v
	=
	\int_{\mathbb R^N}(u^+)^{p-1}v\,dx,
	\qquad
	u,v\in X.
	\]
Thus, after replacing \(\Phi\) by \(\Phi^+\), the preceding
basic nonsmooth variational properties remain valid.
	
We first verify uniform mountain pass geometry for \(I_\lambda^+\).
Since \(0\le \Phi^+(u)\le \Phi(u)\), for every \(\lambda>0\),
	\[
	I_\lambda^+(u)
	=
	\lambda\Psi(u)-\Phi^+(u)
	\ge
	\lambda\Psi(u)-\Phi(u)
	=
	I_\lambda(u).
	\]
Lemma \ref{MP-Q} therefore gives
	\[
	\inf_{\|u\|_X=\rho_0}
	I_{1-\varepsilon_0}^+(u)
	\ge
	\inf_{\|u\|_X=\rho_0}
	I_{1-\varepsilon_0}(u)
	\ge
	\alpha_0>0.
	\]
For the endpoint, choose \(e\ge0\). Then \(\Phi^+(e)=\Phi(e),\)
and hence
	\[
	I_{1+\varepsilon_0}^+(e)
	=
	I_{1+\varepsilon_0}(e)
	<0.
	\]
Thus the families \(I_\lambda^+\) and \(I_\lambda\)
have uniform mountain pass geometry.
	
Writing \(c_+(\lambda)
	=
	\inf_{\gamma\in\Gamma_0}
	\sup_{t\in[0,1]}
	I_\lambda^+(\gamma(t)),\)
we have
	\[
	c_+(\lambda)\ge\alpha_0>0,
	\qquad
	\lambda\in[1-\varepsilon_0,1+\varepsilon_0].
	\]
As in Section 3,
choose a sequence
	\[
	\lambda_n\downarrow1,
	\qquad
	\lambda_n\in[1,1+\varepsilon_0],
	\]
such that \(c_+'(\lambda_n)\) exists. For every \(n\in\mathbb N\),
there is a critical point \(w_n\in D(\Psi),\) of \(I_{\lambda_n}^+\)
satisfying
	\(I_{\lambda_n}^+(w_n)=c_+(\lambda_n).\)
For \(I_\lambda^+\), Lemma \ref{lem:IB-condition}
also applies. Indeed,
	\[
	(u_t)^+(x)
	=
	tu^+\left(\frac{x}{t}\right).
	\]
and critical points of \(I_\lambda^+\) satisfy the identity
	\[
	\lambda
	\left[
	NK(u)+\frac{N+2}{2}B(u)+\frac12D(u)
	\right]
	=
	\frac{N+p}{p}
	\int_{\mathbb R^N}(u^+)^p\,dx.
	\]
This gives \(\sup_{n\in\mathbb N}\|w_n\|_X<+\infty.\) for \(w_n\in D(\Psi)\).

The uniform boundedness statement implies that \(\{w_n\}\) is bounded in \(X\).
Repeating the parameter-limit argument of Lemma \ref{thm:critical-point-lambda-one},
we obtain a subsequence and \(w\in D(\Psi)\) such that \(w_n\to w\) in \(X\).
Since \(w_n\) is critical for \(I_{\lambda_n}^+\),
letting \(n\to\infty\) gives
\begin{equation}\label{i1++}
	\Psi(v)-\Psi(w)-(\Phi^+)'(w)(v-w)\ge0,
	\qquad
	\forall v\in D(\Psi).
\end{equation}
Hence \(w\) is critical for \(I_1^+\).

Write
\[
\bar c_+
=
\lim_{\lambda\downarrow1}c_+(\lambda).
\]
Since \(c_+(\lambda)\ge\alpha_0\),
\[
I_1^+(w)=\bar c_+\ge\alpha_0>0.
\]
Therefore \(w\not\equiv0\).

By Remark \ref{rem:positive-truncation-regularity},
\(w\) weakly solves the truncated equation
	\[
	\operatorname{div}
	\left(
	\frac{Dw}{\sqrt{1-|Dw|^2}}
	\right)
	-
	V(x)w
	+
	(w^+)^{p-1}
	=0
	\]
.
	
We prove \(w\ge0\). Set \(w^-=\max\{-w,0\}.\)
Since \(w\in D(\Psi)\), we have \(w^+\in D(\Psi)\).
Taking \(v=w^+\) in \eqref{i1++} yields
	\[
	\Psi(w^+)-\Psi(w)
	-
	(\Phi^+)'(w)(w^+-w)
	\ge0.
	\]
Observe that \(w^+-w=w^-,\)
and
	\[
	(\Phi^+)'(w)w^-
	=
	\int_{\mathbb R^N}(w^+)^{p-1}w^-\,dx
	=
	0.
	\]
Consequently,
	\[
	\Psi(w^+)-\Psi(w)\ge0.
	\]
On the other hand,
	\[
	\Psi(w^+)-\Psi(w)
	=
	-
	\int_{\{w<0\}}G(Dw)\,dx
	-
	\frac12
	\int_{\{w<0\}}V(x)(w^-)^2\,dx.
	\]
Therefore,
	\[
	-
	\int_{\{w<0\}}G(Dw)\,dx
	-
	\frac12
	\int_{\{w<0\}}V(x)(w^-)^2\,dx
	\ge0.
	\]
Both terms on the left are nonpositive and \(V(x)\ge V_0>0\), so \(w^-=0.\)
Thus \(w\ge0\) in \(\mathbb R^N\).
	
Since \(w\ge0\),
	\[
	(w^+)^{p-1}
	=
	w^{p-1}
	=
	|w|^{p-2}w.
	\]
Consequently \(w\) is a nonnegative weak solution of \eqref{biff}.
	
Section 4 gives \(w\in C_{\mathrm{loc}}^{1,\alpha}\) and \(|Dw|<1\).
On any bounded subdomain, the operator can be written as
	\[
	-\operatorname{div}(a(x)Dw)+V(x)w=w^{p-1}\ge0,
	\qquad
	a(x)=\frac1{\sqrt{1-|Dw|^2}},
	\]
and \(a(x)\) is uniformly bounded and uniformly positive definite there.
The weak Harnack inequality implies that a nonnegative weak solution \(w\) is either identically zero
or strictly positive everywhere. Since \(w\not\equiv0\), we conclude \(w(x)>0, x\in\mathbb R^N.\)
\end{proof}
The symmetric geometry and compactness also give multiplicity.
\begin{theorem}\label{solution-infty}
The original problem \eqref{biff} admits infinitely many distinct nontrivial weak solutions \(\{u^{(k)}\}_{k=1}^{\infty}
	\subset Y.\)
\end{theorem}
\begin{proof}
By Lemma \ref{prop:Phi-Psi-conditions}, \(\Psi\) satisfies \((\Psi_1)\), \((\Psi_2)\), and \((\Psi_3)\),
while \(\Phi\) satisfies \((\Phi_1)\).
Lemma \ref{lem:Phi2-condition} gives \((\Phi_2)\) for \(\Phi\).
Lemma \ref{lem:even-condition} gives evenness of the family \(I_\lambda=\lambda\Psi-\Phi\),
	\((E)\). 
Lemma \ref{lem:SMP-and-divergence} gives uniform symmetric mountain pass geometry (SMP) for \(I_\lambda\).
By Lemma \ref{lem:IB-condition}, the family \(I_\lambda\) also satisfies (IB).
	
Thus \(I_\lambda\) satisfies all the hypotheses of Lemma \ref{th1.7}. Applying Lemma \ref{th1.7},
the original functional \(I_1\) has a sequence of critical points
	\[
	\{u^{(k)}\}_{k=1}^{\infty}
	\subset D(\Psi),
	\]
which can be chosen so that
	\[
	I_1\bigl(u^{(k)}\bigr)=c_k(1),
	\qquad
	c_k(1)\to+\infty, 
	\quad k\to\infty.
	\]
	
Since \(c_k(1)\ge \alpha_0>0,\) and \(I_1(0)=0,\)
we have \(u^{(k)}\not\equiv0.\) for every \(k\in\mathbb N\)
Moreover, \(c_k(1)\to+\infty\)
These critical points cannot take only finitely many values. Otherwise,
\(\{I_1(u^{(k)})\}_{k=1}^{\infty}\) would take only finitely many values, contradicting \(c_k(1)\to+\infty\).
Hence \(\{u^{(k)}\}_{k=1}^{\infty}\) gives infinitely many distinct nontrivial critical points.
	
By Theorem \ref{thm:critical-to-weak-regular}, every critical point of \(I_1\) weakly solves \eqref{biff}.
Thus, for every \(k\in\mathbb N\), \(u^{(k)}\) is a nontrivial weak solution of \eqref{biff}.
	
Since \(u^{(k)}\in D(\Psi)\subset X\subset Y,\)
the original problem \eqref{biff} admits infinitely many distinct nontrivial weak solutions
	\[
	\{u^{(k)}\}_{k=1}^{\infty}\subset Y.
	\]
\end{proof}

\begin{theorem}\label{thm:main-classical-regularity}
The positive weak solution \(w\) of Theorem \ref{solution-positive}
and all weak solutions of Theorem \ref{solution-infty}
\(u^{(k)}\) satisfy, for every \(q>N\),
	\[
	w,u^{(k)}\in W_{\mathrm{loc}}^{2,q}(\mathbb R^N),
	\]
and, in \(\mathbb R^N\),
	\[
	|Dw|<1,
	\qquad
	|Du^{(k)}|<1.
	\]
Moreover, for some \(\alpha\in(0,1)\),
	\[
	w,u^{(k)}\in
	C_{\mathrm{loc}}^{2,\alpha}(\mathbb R^N).
	\]
Thus \(w\) is a positive classical solution of the original equation,
and \(\{u^{(k)}\}_{k=1}^{\infty}\) are infinitely many
nontrivial classical solutions.
\end{theorem}
\begin{proof}
By the proof of Theorem \ref{solution-positive}, \(w\)
is a positive weak solution of \eqref{biff}. By Remark
\ref{rem:positive-truncation-regularity}, for every
\(q>N\), we have \(w\in W_{\mathrm{loc}}^{2,q}(\mathbb R^N),\) and \(|Dw|<1\) in \(\mathbb R^N\).
	
Theorem \ref{solution-infty} also shows that, for every
\(k\in\mathbb N\), \(u^{(k)}\in D(\Psi)\) is a nontrivial critical point of \(I_1\)
. By Theorem
\ref{thm:critical-to-weak-regular}, \(u^{(k)}\)
weakly solves \eqref{biff} and, for every \(q>N\),
\(u^{(k)}\in W_{\mathrm{loc}}^{2,q}(\mathbb R^N),\) and \(|Du^{(k)}|<1\) in \(\mathbb R^N\).
	
We now obtain classical regularity. Let \(U\) be a weak solution of type \(w\) or \(u^{(k)}\).
We have
	\[
	\operatorname{div}\left(\frac{DU}{\sqrt{1-|DU|^2}}\right)
	=
	\sum_{i,j=1}^N a^{ij}(DU)U_{ij},
	\]
where
	\[
	a^{ij}(DU)
	=
	\frac{\delta_{ij}}{\sqrt{1-|DU|^2}}
	+
	\frac{U_i U_j}{(1-|DU|^2)^{3/2}},
	\]
the original equation can be written in nondivergence form as
	\[
	\sum_{i,j=1}^N a^{ij}(DU)U_{ij}
	=
	V(x)U-|U|^{p-2}U.
	\]
For every \(\xi\in\mathbb R^N\),
	\[
	\sum_{i,j=1}^N a^{ij}(DU)\xi_i\xi_j
	=
	\frac{|\xi|^2}{\sqrt{1-|DU|^2}}
	+
	\frac{(DU\cdot\xi)^2}{(1-|DU|^2)^{3/2}}.
	\]
Therefore, on a compact set \(K\),
	\begin{equation}\label{down}
		\sum_{i,j=1}^N a^{ij}(DU)\xi_i\xi_j
		\ge
		|\xi|^2.
	\end{equation}
	
Theorem \ref{thm:critical-to-weak-regular} also gives, for every \(q>N\),
	\[
	U\in W_{\mathrm{loc}}^{2,q}(\mathbb R^N),
	\qquad
	|DU|<1
	\quad
	\text{in }\mathbb R^N.
	\]
We first deduce local continuity of \(DU\). Fix \(R>0\). From
\(U\in W^{2,q}(B_R)\), for every \(i=1,\dots,N\) we have \(	D_i U\in W^{1,q}(B_R).\)
Since \(q>N\), Lemma \ref{lem:morrey-local} gives
	\(D_i U\in C^{0,\alpha}(\overline{B_R}),\)
Thus
	\begin{equation}\label{duc}
		DU\in C_{\mathrm{loc}}^{0,\alpha}(\mathbb R^N;\mathbb R^N),
		\qquad
		U\in C_{\mathrm{loc}}^{1,\alpha}(\mathbb R^N).
	\end{equation}
and \(DU\) is locally continuous.
	
Fix a compact set \(K\Subset\mathbb R^N\). Throughout \(K\),
\(|DU(x)|<1\); continuity in \eqref{duc} and compactness therefore imply
the existence of \(\delta_K>0\)
such that, throughout \(K\),
	\(|DU|\le 1-\delta_K.\)
Thus there is a constant depending only on \(\delta_K\),
\(\Lambda_K>0\), such that
	\begin{equation}\label{up}
		\sum_{i,j=1}^N a^{ij}(DU)\xi_i\xi_j
		\le
		\Lambda_K|\xi|^2.
	\end{equation}
Together, \eqref{down} and \eqref{up} show that the coefficient matrix \((a^{ij}(DU))\) is uniformly elliptic on \(K\).
	
We next verify local H{\"o}lder continuity of the coefficients and source. Since
	\(U\in C_{\mathrm{loc}}^{1,\alpha}(\mathbb R^N),\)
we have
	\(DU\in C_{\mathrm{loc}}^{0,\alpha}(\mathbb R^N)\).
Since the map
	\[
	z\mapsto
	\frac{\delta_{ij}}{\sqrt{1-|z|^2}}
	+
	\frac{z_i z_j}{(1-|z|^2)^{3/2}}
	\]
is smooth on \(\{|z|\le1-\delta_K\}\), we obtain
	\begin{equation}\label{aijholder}
		a^{ij}(DU)\in C^{0,\alpha}(K).
	\end{equation}
The condition \(V\in C^1(\mathbb R^N)\) implies that \(V\) is Lipschitz on \(K\); in particular, \(V\in C^{0,\alpha}(K)\).
Hence
	\(V(x)U\in C^{0,\alpha}(K).\)
Since \(p>2\), the function \(s\mapsto |s|^{p-2}s\)
is locally Lipschitz on bounded intervals, and therefore
	\(|U|^{p-2}U\in C^{0,\alpha}(K).\)
Therefore,
	\begin{equation}\label{vupholder}
		V(x)U-|U|^{p-2}U\in C^{0,\alpha}(K).
	\end{equation}
	
Uniform ellipticity, \eqref{aijholder}, and \eqref{vupholder} allow us to apply local Schauder estimates, yielding
	\(U\in C_{\mathrm{loc}}^{2,\alpha}(\mathbb R^N).\)
In particular, \(U\) satisfies, classically,
	\[
	\operatorname{div}\left(\frac{DU}{\sqrt{1-|DU|^2}}\right)
	-
	V(x)U
	+
	|U|^{p-2}U
	=
	0,
	\qquad x\in\mathbb R^N.
	\]
Since \(U\in Y\) and \(U\not\equiv0\), the function \(U\) is a nontrivial classical solution in the natural space \(Y\).
	
Thus the positive weak solution \(w\) of Theorem \ref{solution-positive} and
each weak solution \(u^{(k)}\) of Theorem \ref{solution-infty} satisfy
	\[
	w,u^{(k)}\in C_{\mathrm{loc}}^{2,\alpha}(\mathbb R^N).
	\]
and hence solve the original problem \eqref{biff} in the classical sense.
\end{proof}

\section*{Acknowledgement}
The research was supported by NSFC-12371107.

\bibliographystyle{plainnat}
\bibliography{references}
\end{document}